\documentclass[12pt]{amsart}
\usepackage{amsfonts,amssymb,color}
\usepackage[mathscr]{eucal}
\usepackage{amsmath, amsthm}
\usepackage{mathrsfs}
\usepackage{amsbsy}
\usepackage{dsfont}
\usepackage{bbm}
\usepackage{wasysym}
\usepackage{stmaryrd}
\usepackage{undertilde}
\usepackage{url}
\input xypic
\xyoption{all}
\begin{document}
\baselineskip = 5mm
\newcommand \ZZ {{\mathbb Z}} 
\newcommand \FF {{\mathbb F}} %
\newcommand \NN {{\mathbb N}} 
\newcommand \QQ {{\mathbb Q}} 
\newcommand \RR {{\mathbb R}} 
\newcommand \CC {{\mathbb C}} 
\newcommand \PR {{\mathbb P}} 
\newcommand \AF {{\mathbb A}} 
\newcommand \bcA {{\mathscr A}}
\newcommand \bcB {{\mathscr B}}
\newcommand \bcC {{\mathscr C}}
\newcommand \bcD {{\mathscr D}}
\newcommand \bcE {{\mathscr E}}
\newcommand \bcF {{\mathscr F}}
\newcommand \bcG {{\mathscr G}}
\newcommand \bcH {{\mathscr H}}
\newcommand \bcM {{\mathscr M}}
\newcommand \bcN {{\mathscr N}}
\newcommand \bcI {{\mathscr I}}
\newcommand \bcJ {{\mathscr J}}
\newcommand \bcK {{\mathscr K}}
\newcommand \bcL {{\mathscr L}}
\newcommand \bcO {{\mathscr O}}
\newcommand \bcP {{\mathscr P}}
\newcommand \bcQ {{\mathscr Q}}
\newcommand \bcR {{\mathscr R}}
\newcommand \bcS {{\mathscr S}}
\newcommand \bcT {{\mathscr T}}
\newcommand \bcU {{\mathscr U}}
\newcommand \bcV {{\mathscr V}}
\newcommand \bcW {{\mathscr W}}
\newcommand \bcX {{\mathscr X}}
\newcommand \bcY {{\mathscr Y}}
\newcommand \bcZ {{\mathscr Z}}
\newcommand \goa {{\mathfrak a}}
\newcommand \gob {{\mathfrak b}}
\newcommand \goc {{\mathfrak c}}
\newcommand \gom {{\mathfrak m}}
\newcommand \gop {{\mathfrak p}}
\newcommand \goT {{\mathfrak T}}
\newcommand \goC {{\mathfrak C}}
\newcommand \goD {{\mathfrak D}}
\newcommand \goM {{\mathfrak M}}
\newcommand \goN {{\mathfrak N}}
\newcommand \goS {{\mathfrak S}}
\newcommand \goH {{\mathfrak H}}
\newcommand \uno {{\mathbbm 1}}
\newcommand \Le {{\mathbbm L}}
\newcommand \Spec {{\rm {Spec}}}
\newcommand \Pic {{\rm {Pic}}}
\newcommand \Jac {{{J}}}
\newcommand \Alb {{\rm {Alb}}}
\newcommand \NS {{{NS}}}
\newcommand \Corr {{Corr}}
\newcommand \Chow {{\mathscr C}}
\newcommand \Sym {{\rm {Sym}}}
\newcommand \Alt {{\rm {Alt}}}
\newcommand \Prym {{\rm {Prym}}}
\newcommand \cone {{\rm {cone}}}
\newcommand \Const {{\rm {Const}}}
\newcommand \cha {{\rm {char}}}
\newcommand \eff {{\rm {eff}}}
\newcommand \tr {{\rm {tr}}}
\newcommand \pr {{\rm {pr}}}
\newcommand \shf {{\rm {a}}}
\newcommand \ev {{\it {ev}}}
\newcommand \Id {{\rm {Id}}}
\newcommand \Nis {{\rm {Nis}}}
\newcommand \et {{\rm {\acute e t}}}
\newcommand \prop {{\rm {prop}}}
\newcommand \equi {{\rm {equi}}}
\newcommand \proeq {{\rm {proeq}}}
\newcommand \CS {{C}}
\newcommand \CMon {\it {CMon}}
\newcommand \AGr {\it {AbGr}}
\newcommand \forget {{\it Forget}}
\newcommand \we {{\rm {w\! .\, e\! .}}}
\newcommand \supp {{\rm Supp}}
\newcommand \interior {{\rm {Int}}}
\newcommand \sep {{\rm {sep}}}
\newcommand \td {{\rm {tdeg}}}
\newcommand \sw {{\rm {sw}}}
\newcommand \alg {{\rm {alg}}}
\newcommand \im {{\rm im}}
\newcommand \Ho {{\sf H}}
\newcommand \dom {{\rm dom}}
\newcommand \codom {{\rm codom}}
\newcommand \holim {{\rm holim}}
\newcommand \Pre {{\mathscr P}}
\newcommand \Funct {{\rm Funct}}
\newcommand \op {{\rm op}}
\newcommand \Hom {{\rm Hom}}
\newcommand \bHom {{\bf {Hom}}}
\newcommand \uHom {{\underline {\rm Hom}}}
\newcommand \cHom {{\mathscr H\! }{\it om}}
\newcommand \map {{\rm {map}}}
\newcommand \Map {{\rm {Map}}}
\newcommand \Hilb {{\rm Hilb}}
\newcommand \Sch {{\mathscr S\! }{\it ch}}
\newcommand \Shv {{\mathscr S\! hv}} 
\newcommand \cHilb {{\mathscr H\! }{\it ilb}}
\newcommand \cExt {{\mathscr E\! }{\it xt}}
\newcommand \colim {{{\rm colim}\, }} 
\newcommand \bSpec {{\bf {Spec}}}
\newcommand \fin {{\rm {f}}}
\newcommand \ind {{\rm {ind}}}
\newcommand \PShv {{\mathscr P\mathscr S\! hv}}
\newcommand \End {{\rm {End}}}
\newcommand \coker {{\rm {coker}}}
\newcommand \coeq {{{\rm coeq}\, }} 
\newcommand \id {{\rm {id}}}
\newcommand \van {{\rm {van}}}
\newcommand \spc {{\rm {sp}}}
\newcommand \Ob {{\rm Ob}}
\newcommand \Aut {{\rm Aut}}
\newcommand \cor {{\rm {cor}}}
\newcommand \res {{\rm {res}}}
\newcommand \Gal {{\rm {Gal}}}
\newcommand \Mon {{\mathscr M\! }{\it on}}
\newcommand \Gr {{\mathscr G\! }{\it r}}
\newcommand \PGL {{\rm {PGL}}}
\newcommand \Bl {{\rm {Bl}}}
\newcommand \Sing {{\it {Sing}}}
\newcommand \spn {{\rm {span}}}
\newcommand \Nm {{\rm {Nm}}}
\newcommand \inv {{\rm {inv}}}
\newcommand \codim {{\rm {codim}}}
\newcommand \ptr {{\pi _2^{\rm tr}}}
\newcommand \sg {{\Sigma }}
\newcommand \Sm {{\mathscr S\! m}} 
\newcommand \CHM {{\mathscr C\! \mathscr M}}
\newcommand \DM {{\sf DM}}
\newcommand \FS {{FS}}
\newcommand \MM {{\mathscr M\! \mathscr M}}
\newcommand \HS {{\mathscr H\! \mathscr S}}
\newcommand \Ex {{\it Ex}}
\newcommand \ExG {{{\it Ex}^{\bcG }}}
\newcommand \RD {{\rm R}}
\newcommand \MHS {{\mathscr M\! \mathscr H\! \mathscr S}}
\newcommand \de {\Delta }
\newcommand \deop {{\Delta \! }^{op}\, }
\newcommand \Sets {{\mathscr S\! ets}}
\newcommand \SSets {{\deop \mathscr S\! ets}}
\newcommand \Vect {{\mathscr V\! ect}}
\newcommand \Gm {{{\mathbb G}_{\rm m}}}
\newcommand \trdeg {{\rm {tr.deg}}}
\newcommand \univ {\tiny {\wasylozenge }}
\newcommand \tame {\rm {tame }}
\newcommand \prym {\tiny {\Bowtie }}
\newcommand \swc {{\Box }}
\newtheorem{theorem}{Theorem}
\newtheorem{lemma}[theorem]{Lemma}
\newtheorem{corollary}[theorem]{Corollary}
\newtheorem{proposition}[theorem]{Proposition}
\newtheorem{remark}[theorem]{Remark}
\newtheorem{definition}[theorem]{Definition}
\newtheorem{conjecture}[theorem]{Conjecture}
\newtheorem{example}[theorem]{Example}
\newtheorem{question}[theorem]{Question}
\newtheorem{assumption}[theorem]{Assumption}
\newtheorem{warning}[theorem]{Warning}
\newtheorem{fact}[theorem]{Fact}
\newtheorem{crucialquestion}[theorem]{Crucial Question}
\newcommand \lra {\longrightarrow}
\newcommand \hra {\hookrightarrow}
\def\blue {\color{blue}}
\def\red {\color{red}}
\def\green {\color{green}}
\newenvironment{pf}{\par\noindent{\em Proof}.}{\hfill\framebox(6,6)
\par\medskip}
\title[$\AF ^1$-connectivity v.s. rational equivalence]
{\bf $\AF ^1$-connectivity on Chow monoids v.s. rational equivalence of algebraic cycles}
\author{Vladimir Guletski\u \i }

\date{07 November 2015}


\begin{abstract}
\noindent Let $k$ be a field of characteristic zero, and let $X$ be a projective variety embedded into a projective space over $k$. For two natural numbers $r$ and $d$ let $C_{r,d}(X)$ be the Chow scheme parametrizing effective cycles of dimension $r$ and degree $d$ on the variety $X$. Choosing an $r$-cycle of minimal degree gives rise to a chain of embeddings of Chow schemes, whose colimit is the connective Chow monoid $C_r^{\infty }(X)$ of $r$-cycles on $X$. Let $BC_r^{\infty }(X)$ be the classifying motivic space of this monoid. In the paper we establish an isomorphism between the Chow group $CH_r(X)_0$ of degree $0$ dimension $r$ algebraic cycles modulo rational equivalence on $X$, and the group of sections of the Nisnevich sheaf of $\AF ^1$-path connected components of the loop space of $BC_r^{\infty }(X)$ at $\Spec (k)$. Equivalently, $CH_r(X)_0$ is isomorphic to the group of sections of the stabilized motivic fundamental group $\Pi _1^{S^1\wedge \AF ^1}(BC_r^{\infty }(X))$ at $\Spec (k)$.
\end{abstract}







\maketitle


\section{Introduction}
\label{intro}

Algebraic cycles are linear combinations of closed integral subschemes in algebraic varieties over a field. Two algebraic cycles $A$ and $B$ on a variety $X$ are said to be rationally equivalent if there exists an algebraic cycle $Z$ on $X\times \PR ^1$, such that, for two fundamental points $0$ and $\infty $ on $\PR ^1$, the cycle-theoretic fibres $Z(0)$ and $Z(\infty )$ coincide with $A$ and $B$ respectively. Rational equivalence is a fundamental notion in algebraic geometry, which substantially depends on the intersection multiplicities tacitly involved into the definition above. Intersection multiplicities are well controlled in cycles which are cascade intersections of cycles starting from codimension one. This is not always the case, of course. For example, if $X$ is a $K3$-surface, the Chow group of $0$-cycles modulo rational equivalence on $X$ is large, in the sense that it cannot be parametrized by an abelian variety over the ground field, \cite{Mumford}. On the other hand, its subgroup generated by divisorial intersections on $X$ is just $\ZZ $, see \cite{BeauvilleVoisin}. This example tells us that intersection multiplicities are geometrically manageable only for a small fraction of all algebraic cycles appearing in nature.

Another difficulty with algebraic cycles is that they are originally given in terms of groups, i.e. positive and negative multiplicities can appear in a cycle simultaneously. The use of negative numbers was questionable for mathematicians dealing with algebraic equations in sixteenth century. In modern terms, the concern can be expressed by saying that the completion of a monoid is a too formal construction. The problem might seem to be not that funny when passing to the completions of Chow monoids, i.e. gatherings of Chow varieties parametrizing effective cycles on projective varieties embedded into projective spaces. The Chow monoids themselves are geometrically given in terms of Cayley forms, whereas their completions are less visible.

These two things taken together have an effect that, in contrast to rational connectivity, rational equivalence is difficult to deform in a smooth projective family over a base, cf. \cite{KMM}. As a consequence, the deep conjectures on rational equivalence are hard to approach, and by now they are solved in a small number of cases (see, for example, \cite{Sur les zero-cycles}). The state of things would be possibly better if we could recode rational equivalence into more effective (i.e. positive) data, appropriate for deformation in smooth projective families over a nice base. The purpose of the present paper is to investigate whether the $\AF ^1$-homotopy type can help in finding such data.

More precisely, let $X$ be a projective variety over a field $k$, and fix an embedding of $X$ into the projective space $\PR ^m$. To avoid the troubles with representability of Chow sheaves in positive characteristic, we must assume that $k$ is of characteristic zero. Effective algebraic cycles of dimension $r$ and degree $d$ on $X$, considered with regard to the embedding $X\subset \PR ^m$, are represented by the Chow scheme $C_{r,d}(X)$ over $k$. Let $Z_0$ be an $r$-cycle of minimal degree $d_0$ on $X$. For example, if $r=0$ then $Z_0$ can be a point, and if $r=1$ then $Z_0$ can be a line on $X$. The cycle $Z_0$ gives rise to a chain of embeddings $C_{r,d}(X)\subset C_{r,d+d_0}(X)$, whose colimit $C_r^{\infty }(X)$ is the {\it connective Chow monoid} of effective $r$-cycles on $X$. Let $C_r^{\infty }(X)^+$ be the group completion of $C_r^{\infty }(X)$ in the category of set-valued simplicial sheaves on the smooth Nisnevich site over $k$. Let also $\Pi _0^{\AF ^1}$ be the functor of $\AF ^1$-connected components and $\Pi _1^{\AF ^1}$ be the functor of the $\AF ^1$-fundamental group on simplicial Nisnevich sheaves, see \cite{MorelVoevodsky} or \cite{AsokMorel}.

Now, consider the Chow group $CH_r(X)_0$ of degree zero $r$-cycles modulo rational equivalence on $X$, where the degree of cycle classes is given with regard to the embedding of $X$ into $\PR ^m$ over $k$. Any finitely generated field extension $K$ of the ground field $k$ is the function field $k(Y)$ of an irreducible variety $Y$ over $k$. For a simplicial sheaf $\bcF $, let $\bcF (K)$ be the stalk of $\bcF $ at the generic point $\Spec (K)$ of the variety $Y$. In the paper we establish a canonical (up to a projective embedding) isomorphism
  $$
  CH_r(X_K)_0\simeq \Pi _0^{\AF ^1}(C_r^{\infty }(X)^+)(K)\; ,
  $$
for an arbitrary finitely generated field extension $K$ over $k$ (Theorem \ref{main1}).

Let, furthermore, $BC_r^{\infty }(X)$ be the motivic classifying space of the connective Chow monoid $C_r^{\infty }(X)$. We also prove that
  $$
  CH_r(X_K)_0\simeq
  \Pi _0^{\AF ^1}(\Omega \Ex BC_r^{\infty }(X))(K)
  \; ,
  $$
where $\Omega $ is right adjoint to the simplicial suspension $\Sigma $ in the pointed category of simplicial Nisnevich sheaves, and $\Ex $ is a fibrant replacement functor for simplicial Nisnevich sheaves (Corollary \ref{main2.5}). Another reformulation of the main result is in terms of $S^1\wedge \AF ^1$-fundamental groups, where $S^1$ is the simplicial circle. Namely,
  $$
  CH_r(X_K)_0\simeq
  \Pi _1^{S^1\wedge \AF ^1}
  (BC_r^{\infty }(X))(K)\; ,
  $$
i.e. the Chow group of $r$-cycles of degree zero modulo rational equivalence on $X$ is isomorphic to the stalk at $\Spec (K)$ of the $S^1\wedge \AF ^1$-fundamental group of the motivic classifying space of the Chow monoid $C_r^{\infty }(X)$ (Corollary \ref{main3}). The smashing by $S^1$ is a sort of stabilization, and not yet fully understood.

The use of the second isomorphism is that it encodes rational equivalence on $r$-cycles in terms of $\AF ^1$-path connectedness on the motivic space $L_{\AF ^1}\Omega \Ex BC_r^{\infty }(X)$. The localization functor $L_{\AF ^1}$ is a transfinite machine, which can be described in terms of sectionwise fibrant replacement, the Godement resolution, homotopy limit of the corresponding cosimplicial simplicial sheaves and the Suslin-Voevodsky's singularization functor. The quadruple operation $L_{\AF ^1}\Omega \Ex B$ is a bigger machine recoding rational equivalence into $\AF ^1$-path connectivity, at some technical cost, of course.

The proof of the main result (Theorem \ref{main1}) is basically a gathering of known facts in $\AF ^1$-homotopy theory of schemes and Chow sheaves, collected in the right way. The substantial arguments are Lemma \ref{Pi0&completion} and the use of Proposition 6.2.6 from the paper \cite{AsokMorel} by Asok and Morel. In Section \ref{concommon} we introduce the needed tools from simplicial sheaves on a small site and the functor $\Pi _0$. Section \ref{localization} is devoted to the Bousfield localization of simplicial sheaves by an interval and to proving Lemma \ref{Pi0&completion} which says that the group completion commutes with the localized $\Pi _0$. In Section \ref{Chowsheaves} we pass to Nisnevich sheaves on schemes and recall the theory of Chow sheaves following \cite{SV-ChowSheaves}. The main results appear in Section \ref{rateq-ratcon}, where we prove the existence of the canonical (up to a projective embedding) isomorphisms between the Chow groups and the stalks of the corresponding motivic homotopy groups of $C_r^{\infty }(X)^+$ and $BC_r^{\infty }(X)$. In Appendix we collect the needed basics from homotopical algebra, in order to make the text more self-contained.

\medskip

{\sc Acknowledgements.} The paper is written in the framework of the EPSRC grant EP/I034017/1. The author is grateful to Aravind Asok for pointing out a drawback in the proof of Corollary \ref{main3} in the first version of the paper, to Sergey Gorchinskiy for explaining how to remove the degree $1$ cycle assumption from the statement of Theorem \ref{main1}, and to both for their interest and useful comments via email and skype. 

\section{$\Pi _0$ and monoids in simplicial sheaves}
\label{concommon}

Let $\Delta $ be the simplex category, i.e. the category whose objects are finite sets $[n]=\{ 0,1,\dots ,n\} $, for all $n\in \NN $, and morphisms $[m]\to [n]$ are order-preserving functions from $[m]$ to $[n]$. Let $\bcS $ be a cartesian monoidal category with a terminal object $*$. The category $\deop \bcS $ of simplicial objects in $\bcS $ is the category of contravariant functors from $\Delta $ to $\bcS $. Since $[0]$ is the terminal object in $\Delta $, the functor $\Gamma :\deop \bcS \to \bcS $, sending $\bcX $ to $\bcX _0$, is the functor of global sections on simplicial objects in $\bcS $ considered as presheaves on $\Delta $. The functor $\Gamma $ admits left adjoint $\Const :\bcS \to \deop \bcS $ sending an object $\bcX $ in $\bcS $ to the constant presheaf on $\Delta $ determined by $\bcX $.

Assume, moreover, that $\bcS $ is cocomplete. For any object $\bcX $ in $\deop \bcS $, let $\Pi _0(\bcX )$ be the coequalizer of the morphisms $\bcX _1 \rightrightarrows \bcX _0$ induced by the two morphisms from $\Delta [0]$ to $\Delta [1]$. This gives a functor $\Pi _0:\deop \bcS \to \bcS $ and the canonical epimorphism $\Psi :\Gamma \to \Pi _0$. If $\bcY $ is an object in $\bcS $, and $f:\bcX \to \Const (\bcY )$ is a morphism in $\deop \bcS $, the precompositions of $f_0:\bcX _0\to \bcY $ with the two morphisms from $\bcX _1$ to $\bcX _0$ coincide. By universality of the coequalizer, we obtain the morphism $f':\Pi _0(\bcX )\to \bcY $. The correspondence $f\mapsto f'$ is one-to-one and natural in $\bcX $ and $\bcY $. In other words, $\Pi _0$ is left adjoint to $\Const $. Since products in $\deop \bcS $ are objectwise, the functor $\Pi _0$ preserves finite products. Certainly, if $\bcC $ is the terminal category, then $\Pi _0$ is the usual functor of connected components on simplicial sets.

Let $\bcC $ be an essentially small category and let $\tau $ be a subcanonical topology on it. Assume also that $\bcC $ contains all finite products and let $*$ be the terminal object in it. Let $\bcP $ be the category of presheaves of sets on $\bcC $ and let $\bcS $ be the full subcategory of set valued sheaves on $\bcC $ in the topology $\tau $. Since $\tau $ is subcanonical, the Yoneda embedding $h:\bcC \to \bcP $, sending an object $X$ to the representable presheaf $h_X=\Hom _{\bcC }(-,X)$, and a morphism $f:X\to Y$ to the morphism of presheaves $h_f=\Hom _{\bcC }(-,f)$, takes its values in the category of sheaves $\bcS $. If $*$ is the terminal object in $\bcC $, then $h_*=\Hom _{\bcC }(-,*)$ is the terminal object in $\bcP $ and $\bcS $. Limits in $\bcS $ are limits in $\bcP $. In particular, we have objectwise finite products in $\bcS $ and the category $\bcS $ is Cartesian monoidal.

For a presheaf $\bcX $, let $\bcX ^{\shf }$ be the sheaf associated to $\bcX $ in $\tau $. Since $\bcP $ is complete, the sheafification of colimits in $\bcP $ shows that $\bcS $ is cocomplete too. In order to make a difference between $\Pi _0$ in $\deop \bcS $ and $\Pi _0$ in $\deop \bcP $, we shall denote the latter functor by $\pi _0$, so that, for a simplicial sheaf $\bcX $, one has $\Pi _0(\bcX )=\pi _0(\bcX )^{\shf }$. As the coequalizer $\pi _0$ is sectionwise, $\Pi _0(\bcX )$ is the sheafififcation of the presheaf sending $U$ to $\pi _0(\bcX (U))$.

Let $\SSets $ be the category of simplicial sets. For a natural number $n$ let $\Delta [n]$ be the representable functor $\Hom _{\Delta }(-,[n])$. For any sheaf $\bcF $ on $\bcC $ let $\Delta _{\bcF }[n]$ be the simplicial sheaf defined by the formula
  $$
  (\Delta _{\bcF }[n])_m(U)=
  \bcF (U)\times \Hom _{\Delta }([m],[n])\; ,
  $$
for any $U\in \Ob (\bcC )$ and any natural number $m$. This gives the full and faithful embeddings $\Delta _{?}[n]:\bcS \to \deop \bcS $ and $\Delta _{\bcF }[?]:\Delta \to \deop \bcS $. If $\bcF $ is $h_X$, for some object $X$ in $\bcC $, then we write $\Delta _X[n]$ instead of $\Delta _{\bcF }[n]$, and use $\Delta [n]$ instead of $\Delta _*[n]$. To simplify notation further, we shall identify $\bcC $ with its image in $\deop \bcS $ under the embedding $\Delta _{?}[0]$. For example, for any object $X$ in $\bcC $ it is the same as the corresponding simplicial sheaf $\bcX =\Delta _X[0]=\Const (h_X)$, and the same on morphisms in $\bcC $. The cosimplicial object $\Delta [?]:\Delta \to \deop \bcS $ determines the embedding of simplicial sets into $\deop \bcS $, so that we may also identify $\SSets $ with its image in $\deop \bcS $. This gives the structure of a simplicial category on $\deop \bcS $, such that, for any two simplicial sheaves $\bcX $ and $\bcY $,
  $$
  \bHom (\bcX ,\bcY )=
  \Hom _{\deop \bcS }(\bcX \times \Delta [?],\bcY )\; .
  $$
The corresponding (right) action of $\SSets $ on $\deop \bcS $ is given by the formula
  $$
  (\bcX \times K)_n(U)=\bcX _n(U)\times K_n\; ,
  $$
for any simplicial sheaf $\bcX $ and simplicial set $K$. For simplicity of notation, we shall write $\Delta [n]$ instead of $\Delta [n]$. Then $\Delta _X[n]$ is the product of $\Delta _X[0]$ and $\Delta [n]$.

Looking at $\deop \bcS $ as a symmetric monoidal category with regard to the categorical product in it, one sees that it is closed symmetric monoidal. The internal Hom, bringing right adjoint to the Cartesian products, is given by the formula
  $$
  \cHom (\bcX ,\bcY )_n(U)=
  \Hom _{\deop \bcS }(\bcX \times \Delta _U[n],\bcY )\; .
  $$

Throughout the paper we will be working with monoids in $\deop \bcS $. All monoids and groups will be commutative by default. If $\bcX $ is a monoid in $\deop \bcS $, let $\bcX ^+$ be the group completion of $\bcX $ in $\deop \bcS $. The terms of $\bcX ^+$ are the sheaves associated with the sectionwise completions of the terms of $\bcX $ in $\bcP $. One has a morphism from $\bcX \times \bcX $ to $\bcX ^+$, which is an epimorphism in $\deop \bcS $. When no confusion is possible, the termwise and section-wise completion of $\bcX $ in $\deop \bcP $ will be denoted by the same symbol $\bcX ^+$.

Monoids form a subcategory in $\bcP $. The corresponding forgetful functor has left adjoint sending presheaves to free monoids with concatenation as monoidal operation. The notion of a cancellation monoid in $\bcP $ is standard and sectionwise. A free monoid in $\bcP $ is a cancellation monoid. As limits and colimits in $\deop \bcS $ are termwise, the functors $\Gamma $ and $\Const $ preserve monoids and groups and $\Gamma (\bcX ^+)$ is the same as $\Gamma (\bcX )^+$. Since $\Pi _0$ commutes with finite products, it follows that $\Pi _0$ also preserves monoids and groups.

The monoid of natural numbers $\NN $ is a simplicial sheaf on $\bcC $. A pointed monoid in $\deop \bcS $ is a pair $(\bcX ,\iota )$, where $\bcX $ is a monoid in $\deop \bcS $ and $\iota $ is a morphism of monoids from $\NN $ to $\bcX $. A graded pointed monoid is a triple $(\bcX ,\iota ,\sigma )$, where $(\bcX ,\iota )$ is a pointed monoid and $\sigma $ is a morphism of monoids from $\bcX $ to $\NN $, such that $\sigma \circ \iota =\id _{\NN }$, see page 126 in \cite{MorelVoevodsky}. Notice that to define a morphism from $\NN $ to $\bcX $ is equivalent to choose an element in $\bcX _0(*)$.

Let $(\bcX ,\iota ,\sigma )$ be a pointed graded monoid in $\deop \bcS $. Since $\sigma \circ \iota =\id _{\NN }$, it follows that, for any natural $n$ and any object $U$ in $\bcC $, we have two maps $\iota _{U,n}:\NN \to \bcX _n(U)$ and $\sigma _{U,n}:\bcX _n(U)\to \NN $. It implies that $\bcX _n(U)$ is the coproduct of the sets $\sigma _{U,\, n}^{-1}(d)$, for all $d\geq 0$. The sets $\sigma _{U,\, n}^{-1}(d)$ give rise to the simplicial sheaf which we denote by $\bcX ^d$. Then $\bcX $ is the coproduct of $\bcX ^d$ for all $d\geq 0$. The addition of $\iota (1)$ in $\bcX $ induces morphisms of simplicial sheaves $\bcX ^d\to \bcX ^{d+1}$ for all $d\geq 0$. Let $\bcX ^{\infty }$ be the colimit
  $$
  \bcX ^{\infty }=
  \colim (\bcX ^0\to \bcX ^1\to \bcX ^2\to \dots )
  $$
in $\deop \bcS $. Equivalently, $\bcX ^{\infty }$ is the coequalizer of the addition of $\iota (1)$ in $\bcX $ and the identity automorphism of $\bcX $.

Since now we shall assume that the topos $\bcS $ has enough points, and the category $\bcC $ is Noetherian. Since filtered colimits commute with finite products, $\bcX ^{\infty }$ is the colimit taken in the category of simplicial presheaves, i.e. there is no need to take its sheafification. The commutativity of filtered colimits with finite products also yields the canonical isomorphism between the colimit of the obvious diagram composed by the objects $\bcX ^d\times \bcX ^{d'}$, for all $d,d'\geq 0$, and the product $\bcX ^{\infty }\times \bcX ^{\infty }$. Since the colimit of that diagram is the colimit of its diagonal, this gives the canonical morphism from $\bcX ^{\infty }\times \bcX ^{\infty }$ to $\bcX ^{\infty }$. The latter defines the structure of a monoid on $\bcX ^{\infty }$, such that the canonical morphism
  $$
  \pi :\bcX =\coprod _{d\geq 0}\bcX ^d\to \bcX ^{\infty }
  $$
is a homomorphism of monoids in $\deop \bcS $. We call $\bcX ^{\infty }$ the {\it connective} monoid associated to the pointed graded monoid $\bcX $.

Notice that the category of simplicial sheaves is exhaustive. In particular, if all the morphisms $\bcX ^d\to \bcX ^{d+1}$ are monomorphisms, the transfinite compositions $\bcX ^d\to \bcX ^{\infty }$ are monomorphisms too. This happens if $\bcX $ is a termwise sectionwise cancelation monoid, in which case $\bcX ^{\infty }$ is a termwise sectionwise cancelation monoid too.

The above homomorphisms $\pi $ and $\sigma $ give the homomorphism $(\pi ,\sigma )$ from $\bcX $ to $\bcX ^{\infty }\times \NN $. Passing to completions we obtain the homomorphism $(\pi ^+,\sigma ^+)$ from $\bcX ^+$ to $(\bcX ^{\infty })^+\times \ZZ $.

\begin{lemma}
\label{goraimysh}
Assume $\bcX $ is a sectionwise cancelation monoid. Then
  $$
  (\pi ^+,\sigma ^+):\bcX ^+\to (\bcX ^{\infty })^+\times \ZZ
  $$
is an isomorphism.
\end{lemma}

\begin{pf}
Since the site $\bcS $ has enough points, it suffices to prove the lemma sectionwise and termwise. Then, without loss of generality, we may assume that $\bcX $ is a set-theoretical pointed graded cancelation monoid. Clearly, $\iota ^+$ is an injection, $\pi ^+$ is a surjection, and $\pi ^+\iota ^+=0$. Since $\bcX $ is a cancelation monoid, $\bcX ^+$ is the quotient-set of the set $\bcX \times \bcX $ modulo an equivalence relation
  $$
  (x_1,x_2)\sim (x_1',x_2')\Leftrightarrow
  x_1+x_2'=x_2+x_1'\; .
  $$
For any element $(x_1,x_2)$ in $\bcX \times \bcX $ let $[x_1,x_2]$ be the corresponding equivalence class. Since $\bcX $ is a cancelation monoid, so is the monoid $\bcX ^{\infty }$ too. If $\pi ^+[x_1,x_2]$ is zero, that is $[\pi (x_1),\pi (x_2)]=[0,0]$ in $(\bcX ^{\infty })^+$, it is equivalent to say that $\pi (x_1)=\pi (x_2)$. The latter equality means that there exists a positive integer $n$, such that $x_2=x_1+n\iota (1)$, i.e. $[x_1,x_2]=[0,n\iota (1)]$ in $\bcX ^+$. The element $[0,n\iota (1)]$ sits in the image of $\iota ^+$.
\end{pf}

\section{Homotopy completion and localization of $\Pi _0$}
\label{localization}

All the above considerations were categorical. Let us now switch to homotopical algebra and consider the injective model structures on $\deop \bcS $. Recall that a point $P$ of a topos $\bcT $ is an adjoint pair of functors, $P^*:\bcT \to \Sets $ and $P_*:\Sets \to \bcT $, such that $P^*$ is left adjoint to $P_*$ and preserves finite limits in $\bcT $. If $\bcX $ is an object of $\bcT $, then $\bcX _P=P^*(X)$ is the stalk of $\bcX $ at the point $P$. We will assume that the topos $\bcS $ has enough points. Recall that it means that there exists a set of points $P(\bcS )$ of the topos $\bcS $, such that a morphism $f:\bcX \to \bcY $ in $\bcS $ is an isomorphism in $\bcS $ if and only if, for any point $P\in P(\bcS )$, the morphism $f_P:\bcX _P\to \bcY _P$, induced on stalks, is an isomorphism in the category $\Sets $. Respectively, a morphism $f:\bcX \to \bcY $ in $\deop \bcS $ is an isomorphism in $\deop \bcS $ if and only if, for each $P\in P(\bcS )$, the morphism $f_P:\bcX _P\to \bcY _P$ is an isomorphism in $\deop \Sets $.

Now, a morphism $f:\bcX \to \bcY $ in the category of simplicial sheaves $\deop \bcS $ is a weak equivalence in $\deop \bcS $ if and only if for any point $P^*:\bcS \to \Sets $ of the topos $\bcS $ the induced morphism $\deop P^*(f)$ on stalks is a weak equivalence of simplicial sets. Cofibrations are monomorphisms, and fibrations are defined by the right lifting property in the standard way, see Definition 1.2 on page 48 in \cite{MorelVoevodsky}. The pair $(\bcS ,\bcM )$ is then a model category of simplicial sheaves on $\bcC $ in $\tau $. Notice that the model structure $\bcM $ is left proper, see Remark 1.5 on page 49 in loc.cit. One can also show that it is cellular. Let $\Ho $ be the homotopy category $Ho(\deop \bcS )$ of the category $\deop \bcS $ with regard to $\bcM $. For any two simplicial sheaves $\bcX $ and $\bcY $ the set of morphisms from $\bcX $ to $\bcY $ in $\Ho $ will be denoted by $[\bcX ,\bcY ]$.

The simplicial structure on $\deop \bcS $ is compatible with the model one, so that $\bcS $ is a simplicial model category. Since $$
  [\bcX ,\bcY ]\simeq \pi _0\bHom (\bcX ,\bcY )
  $$
and
  $$
  \bHom (\Delta _U[0],\bcX )\simeq \bcX (U)\; ,
  $$
$\Pi _0(\bcX )$ is the sheafififcation of the presheaf
  $$
  \pi _0(\bcX ):U\mapsto \pi _0\bHom (\Delta _U[0],\bcX )=
  [\Delta _U[0],\bcX ]=[\Const (h_U),\bcX ]
  $$
on $\bcC $ in the topology $\tau $. The multiplication of simplicial sheaves and their morphisms by a simplicial set admits right adjoint, so that it commutes with colimits. In particular, $\Pi _0(\Delta _X[n])\simeq \Delta _X[0]$.

A pointed simplicial sheaf $(\bcX ,x)$ is a pair consisting of a simplicial sheaf $\bcX $ and a morphism $x$ from $*$ to $\bcX $. The definition of a morphism of pointed simplicial sheaves is obvious. Let $\deop \bcS _*$ be the category of pointed simplicial sheaves. The corresponding forgetful functor has the standard left adjoint sending $\bcX $ to the coproduct $\bcX _+$ of $\bcX $ and $*$. The model structure $\bcM $ induces the corresponding model structure on $\deop \bcS _*$, such that the above adjunction is a Quillen adjunction. Having two pointed simplicial sheaves $(\bcX ,x)$ and $(\bcY ,y)$, their wedge product $(\bcX ,x)\vee (\bcY ,y)$ is the coproduct in $\deop \bcS _*$, and the smash product $(\bcX ,x)\wedge (\bcY ,y)$ is the contraction of the wedge product in $(\bcX \times \bcY ,x\times y)$.

Let now $S^1$ be the simplical circle $\Delta [1]/\partial \Delta [1]$ pointed by the image of the boundary $\partial \Delta [1]$ in then quotient simplicial set, and let $S^1$ be its image in $\deop \bcS _*$. Define the simplicial suspension endofunctor $\Sigma $ on $\deop \bcS _*$ sending $(\bcX ,x)$ to $S^1\wedge (\bcX ,x)$. Its left adjoint is the simplicial loop functor $\Omega $ sending $(\bcX ,x)$ to $\cHom _*(S^1,(\bcX ,x))$, where $\cHom _*(-,-)$ is the obvious internal Hom in $\deop \bcS _*$.

Let $\bcX $ be a monoid in $\deop \bcS $. For any object $U$ in $\bcC $ and any positive integer $n$ let $N(\bcX _n(U))$ be the nerve of $\bcX _n(U)$, and let $B\bcX $ be the diagonal of the bisimplicial sheaf $\deop \times \deop \to \bcS $ sending $[m]\times [n]$ to the sheaf $U\mapsto N(\bcX _n(U))_m$. Then $(B\bcX )_n$ is $\bcX _n^{\times n}$ for $n>0$ and, by convention, $(B\bcX )_0$ is the terminal object $*$ in $\bcS $, see page 123 in \cite{MorelVoevodsky}. If $\bcC $ is a terminal category, then $B\bcX $ is the usual classifying space of a simplicial monoid $\bcX $ (that is, a monoid in the category of simplicial sets $\deop \Sets $). Just as in topology, there exists a canonical morphism from $\bcX $ to $\Omega B(\bcX )$, which is a weak equivalence if $\bcX $ is a group, loc.cit.

Following Quillen, \cite{QuillenGroupCompletion}, we will say that a simplicial monoid $\bcX $ is good if the morphism $B\bcX \to B\bcX ^+$, induced by the canonical morphism from $\bcX $ to $\bcX ^+$, is a weak equivalence in $\deop \Sets $. If $X$ is a set-theoretical monoid, then $X$ is good if the corresponding constant simplicial monoid $\bcX =\Const (X)$ is good as a simplicial monoid. If $X$ is a free monoid in $\Sets $, then $\bcX =\Const (X)$ is good in $\deop \Sets $, see Proposition Q.1 in loc.cit.

Recall that, for any point $P$ of the topos $\bcS $, the functor $P^*:\bcS \to \Sets $ preserves finite limits. It follows that, if $\bcX $ is a simplicial sheaf monoid, then the stalk $(B\bcX )_P$ of the classifying space $B\bcX $ at $P$ is canonically isomorphic to the classifying space $B(\bcX _P)$ of the stalk $\bcX _P$ of the simplicial sheaf $\bcX $ at $P$. We will say that a simplicial sheaf monoid $\bcX $ is {\it pointwise good}, if the morphism $(B\bcX )_P\to (B\bcX ^+)_P$, is a weak equivalence of simplicial sets for each point $P$ in $P(\bcS )$. This is, of course, equivalent to saying that the morphism $B\bcX \to B\bcX ^+$ is a weak equivalence in $\deop \bcS $, with regard to the model structure $\bcM $.

Now, if $\bcX $ is a monoid in $\bcS $, we will say that $\bcX $ is {\it pointwise free} if $\bcX _P$ is a free monoid in $\Sets $ for each point $P$ in $P(\bcS )$. If $\bcX $ is pointwise free, it does not necessarily mean that $\bcX $ is a free monoid in the category $\bcS $. It is important, however, that if $\bcX _0$ is a pointwise free monoid in $\bcS $, the corresponding constant simplicial sheaf monoid $\bcX =\Const (\bcX _0)$ is pointwise good, which is a straightforward consequence of the first part of Proposition Q.1 in \cite{QuillenGroupCompletion}.

Similarly, we will say that a monoid $\bcX $ in $\bcS $ is a {\it pointwise cancellation} monoid if $\bcX _P$ is a cancellation monoid in $\Sets $ for each point $P$ in $P(\bcS )$. If $\bcX _0$ is a pointwise cancellation monoid, then the simplicial sheaf monoid $\bcX =\Const (\bcX _0)$ is pointwise good by the second part of Quillen's proposition above.

Let $\Ex $ be the fibrant replacement functor $\ExG $, for the model structure $\bcM $, constructed by taking the composition of the sectionwise fibrant replacement of simplicial sets, the Godement resolution and the homotopy limit of the corresponding cosimplicial simplicial sheaf, as in Section 2.1 in \cite{MorelVoevodsky}. Since $\Ex $ preserves finite limits, it preserves monoids and groups. For the same reason, $\Ex $ commutes with taking the classifying spaces of monoids and groups. The right derived functor of $\Omega $ can be computed by precomposing it with $\Ex $. We will need the following variation of Lemma 1.2 on page 123 in \cite{MorelVoevodsky}.

\begin{lemma}
\label{keylemma}
If $\bcX $ is pointwise good, there is a canonical isomorphism
  $$
  \bcX ^+\simeq \Omega \Ex B(\bcX )
  $$
in the homotopy category $\Ho $.
\end{lemma}

\begin{pf}
Since $\bcX $ is pointwise good, the morphism from $B\bcX $ to $B\bcX ^+$, induced by the canonical morphism from $\bcX $ to $\bcX ^+$, is a weak equivalence in $\deop \bcS $. Applying the right derived functor $\RD \Omega $ to the weak equivalence $B\bcX \to B\bcX ^+$ and reverting the corresponding isomorphism in $\Ho $, we obtain the canonical isomorphism from $\Omega \Ex B(\bcX ^+)$ to $\Omega \Ex B(\bcX )$, in the homotopy category $\Ho $. The composition of the canonical morphism $\Sigma \bcX ^+\to B\bcX ^+$ with the weak equivalence $B\bcX ^+\to \Ex B\bcX ^+$ corresponds to the morphism $\bcX ^+\to \Omega \Ex B\bcX ^+$ under the adjunction between $\Sigma $ and $\Omega $. The latter morphism is the composition of the canonical morphism $\bcX ^+\to \Omega B\bcX ^+$ and the morphism $\Omega B\bcX ^+\to \Omega \Ex B\bcX ^+$. The morphism $\bcX ^+\to \Omega B\bcX ^+$ is a weak equivalence because $\bcX ^+$ is a group. Since any simplicial sheaf of groups $\bcG $ can be replaced, up to a weak equivalence, by a fibrant simplicial sheaf of groups, without loss of generality we may assume that $\bcX ^+$ is fibrant (see, for example, Lemma 2.32 on page 83 in \cite{MorelVoevodsky}). Replacing the functor $B$ by the universal cocycle construction $\overline W$, we see that $B$ preserves, up to a weak equivalence, fibrant objects by Theorem 31 in \cite{JardineFieldsLectures}. Then the morphism $\Omega B\bcX ^+\to \Omega \Ex B\bcX ^+$ is a weak equivalence too. Thus, we obtain an isomorphism from $\bcX ^+$ to $\Omega \Ex B\bcX ^+$ in $\Ho $.
\end{pf}

Next, let $A$ be an object of $\bcC $, and let $\bcA $ be the corresponding constant simplicial sheaf $\Delta _A[0]=\Const (h_A)$ in $\deop \bcS $. As in Appendix below, let
  $$
  S=\{ \bcX \wedge \bcA \to \bcX \; |\; \bcX \in \dom(I)\cup \codom (I)\}
  $$
be the set of morphisms induced by the morphism from $\bcA $ to $*$, where $\dom(I)$ and $\codom (I)$ are the sets of domains and codomains of the generating cofibrations in $\bcM $ on $\bcS $. As $\deop \bcS $ is left proper simplicial cellular model category, there exists the left Bousfield localization of $\bcM $ by $S$ in the sense of Hirschhorne, see \cite{Hirsch}. Denote the localized model structure by $\bcM _A$, and let $L_A$ be the corresponding $S$-localization functor, which is a fibrant approximation in $\bcM _A$ on $\deop \bcS $, see Section 4.3 in \cite{Hirsch}, and the earlier work \cite{GoerssJardineLoc}. Let
  $$
  l:\Id _{\deop \bcS }\to L_A
  $$
be the corresponding natural transformation. For any simplicial sheaf $\bcX $ the morphism $l_{\bcX }:\bcX \to L_A(\bcX )$ is a weak cofibration and $L_A(\bcX )$ is $A$-local, i.e. fibrant in $\bcM _A$. The basics on localization functors see Section 4.3 in Hirschhorn's book \cite{Hirsch} and Appendix below.

Let $\Ho _A$ be the homotopy category of simplicial sheaves converting weak equivalences in $\bcM _A$ into isomorphisms. As simplicial sheaves with respect to $\bcM $ form a simplicial closed cartesian monoidal model category, so is the category of simplicial sheaves with respect to $\bcM _A$. All simplicial sheaves are cofibrant, in $\bcM $ and in $\bcM _A$. It follows that the canonical functors from simplicial sheaves to $\Ho $ and $\Ho _A$ are monoidal. See Appendix for more details on all such things. For any two simplicial sheaves $\bcX $ and $\bcY $ let $[\bcX ,\bcY ]_A$ be the set of morphisms from $\bcX $ to $\bcY $ in $\Ho _A$.

Recall that an object $I$ of a category $\bcD $ with a terminal object $*$ is called an interval if there exists a morphism
  $$
  \mu :I\wedge I\to I
  $$
and two morphisms $i_0,i_1:*\rightrightarrows I$, such that
  $$
  \mu \circ (\id _I\wedge i_0)=i_0\circ p
  \quad \hbox{and}\quad
  \mu \circ (\id _I\wedge i_1)=\id _I\; ,
  $$
where $p$ is the unique morphism from $I$ to $*$, and $i_0\coprod i_1:*\coprod *\to I$ is a monomorphism in $\bcD $, see \cite{MorelVoevodsky}. Certainly, the object $A$ is an interval in $\bcC $ if and only if the object $\bcA $ is an interval in $\deop \bcS $.

Since now we shall assume that $A$ is an interval in $\bcC $. The monoidal multiplication by $\bcA $ is a natural cylinder functor on $\deop \bcS $. If $f,g:\bcX \rightrightarrows \bcY $ are two morphisms from $\bcX $ to $\bcY $ in $\deop \bcS $, a left $A$-homotopy from $f$ to $g$ is a morphism $H:\bcX \times \bcA \to \bcY $, such that $H\circ (\id _{\bcX }\times i_0)=f$ and $H\circ (\id _{\bcX }\times i_1)=g$. Since all simplicial sheaves are cofibrant in both model structures $\bcM $ and $\bcM _A$, $A$-homotopy is an equivalence relation on the set $\Hom _{\deop \bcS }(\bcX ,\bcY )$, see \cite{Hovey}, Proposition 1.2.5 (iii). Let $\Hom _{\deop \bcS }(\bcX ,\bcY )_A$ be the set of equivalence classes modulo this equivalence relation. Whenever $\bcY $ is $A$-local, the set $[\bcX ,\bcY ]_A$ is in the natural bijection with the set $\Hom _{\deop \bcS }(\bcX ,\bcY )_A$.

A point of a simplicial sheaf $\bcX $ is, by definition, a morphism from the terminal simplicial sheaf to $\bcX $. Such morphisms can be identified with the set $\bcX _0(*)$. Two points on $\bcX $ are said to be $A$-path connected if and only if they are left homotopic with respect to $\bcA $.

Since $A$ is an interval, the $A$-localizing functor $L_A$ can be chosen to be more explicit than the construction given in \cite{Hirsch}. Following \cite{MorelVoevodsky}, see page 88, we consider the cosimplicial sheaf
   $$
   \Delta _{A^{\bullet }}[0]:\Delta \to \bcS
   $$
sending $[n]$ to the $n$-product
  $$
  (\Delta _A[0])^n=\Delta _{A^n}[0]
  $$
and acting on morphisms as follows. For any morphism $f:[m]\to [n]$ define a morphism of sets
  $$
  f':\{ 1,\dots ,n\} \to \{ 0,1,\dots ,m+1\}
  $$
setting
  $$
  f'(i)=\left\{
  \begin{array}{ll}
  \min \{ l\in \{ 0,\dots ,m\} \; |\; f(l)\geq i\} \; , & \mbox{if this set is nonempty} \\
  m+1 & \mbox{otherwise. }
  \end{array}
  \right.
  $$
If now $\pr _k:A^n\to A$ is the $k$-th projection and $p:A^n\to *$ the unique morphism to the terminal object, where $A^n$ is the $n$-fold product of $A$, then
  $$
  \pr _k\circ \Delta _{A^{\bullet }}[0](f)=\left\{
  \begin{array}{ll}
  \pr _{f'(k)}\; , & \mbox{if}\; f'(k)\in \{ 1,\dots ,m\} \\
  i_0\circ p \; , & \mbox{if}\; f'(k)=m+1 \\
  i_1\circ p\; , & \mbox{if}\; f'(k)=0\; .
  \end{array}
  \right.
  $$
For any $\bcX $ let $\Sing _A(\bcX )$ be the Suslin-Voevodsky simplicial sheaf
  $$
  [n]\mapsto \cHom (\Delta _{A^n}[0],\Delta _{\bcX _n}[0])\; ,
  $$
where the internal $\cHom $ is taken in the category of sheaves $\bcS $. It is functorial in $\bcX $ and $p:A^n\to *$ induces the morphism
  $$
  s:\Id _{\deop \bcS }\to \Sing _A\; .
  $$

Notice that, although the virtue of $A$ to be an interval is not explicitly used in the Suslin-Voevodsky's construction above, it is used in proving the numerous nice properties of the functor $\Sing _A$, see \cite{MorelVoevodsky}. In particular, each morphism $s_{\bcX }$ from $\bcX $ to $\Sing _A(\bcX )$ is an $A$-local weak equivalence, i.e. a weak equivalence with regard to the model structure $\bcM _A$, see Corollary 3.8 on page 89 in loc.cit.

As it is shown in \cite{MorelVoevodsky}, there exists a sufficiently large ordinal $\omega $, such that $L_A$ can be taken to be the composition
  $$
  L_A=(\Ex \circ \Sing _A)^{\omega }\circ \Ex \; ,
  $$
where $\Ex $ is the functor $\ExG $, i.e. the composition of the sectionwise fibrant replacement, the Godement resolution and the homotopy limit of the corresponding cosimplicial simplicial sheaf (see above). Such constructed localization functor $L_A$ is quite explicit, which gives a clearer picture of what are the functors $\pi _0^A$ and $\Pi _0^A$.

The canonical functor from $\deop \bcS $ to $\Ho $ preserves products. In other words, if $\bcX \times \bcY $ is the product of two simplicial sheaves, the same object $\bcX \times \bcY $, with the homotopy classes of the same projections, is the product of $\bcX $ and $\bcY $ in $\Ho $ and in $\Ho _A$ (see Appendix). The advantage of the above explicit $L_A$ is that it commutes with finite products, see Theorem 1.66 on pages 69 - 70 and the remark on page 87 in \cite{MorelVoevodsky}. Most likely, the general Hirschhorne's construction (see Section 4.3 in \cite{Hirsch}) also enjoys this property, but we could not find the proof in the literature.

\begin{remark}
\label{filin}
{\rm The left derived to any localization functor $L_A$ from $\deop \bcS $ to $A$-local objects in $\deop \bcS $ is left adjoint to the right derived of the forgetful functor in the opposite direction on the homotopy level, see Theorem 2.5 on page 71 in \cite{MorelVoevodsky}. This implies, in particular, that any two localizations $L_A$ and $L'_A$ are weak equivalent to each other. Therefore, in all considerations up to (pre-$A$-localized) weak equivalence in $\deop \bcS $ we may freely exchange the localization functor $L_A$ considered in \cite{Hirsch} by the concrete Suslin-Voevodsky's one, and vice versa.
}
\end{remark}

\begin{lemma}
\label{triton}
For any simplicial sheaf $\bcX $ the canonical map
  $$
  \Hom _{\deop \bcS }(*,\bcX )_A
  \to
  \Hom _{\deop \bcS }(*,L_A(\bcX ))_A
  $$
is surjective.
\end{lemma}

\begin{pf}
We know that the natural transformation $l:\Id \to L_A$ induces the epimorphism $\Pi _0\to \Pi _0^A$ by Corollary 3.22 in \cite{MorelVoevodsky}. The morphism $\Psi :\Gamma \to \Pi _0$ is an epimorphism too. This gives that the map
  $$
  \Hom _{\bcS }(*,\bcX _0)
  \to
  \Hom _{\bcS }(*,\Pi _0(\bcX ))
  \to
  \Hom _{\bcS }(*,\Pi _0^A(\bcX ))
  $$
is surjective. By adjunction, $\Hom _{\bcS }(*,\bcX _0)\simeq \Hom _{\deop \bcS }(*,\bcX )$, and since $L_A(\bcX )$ is $A$-local, $\Hom _{\bcS }(*,\Pi _0^A(\bcX ))$ is isomorphic to $\Hom _{\deop \bcS }(*,L_A(\bcX ))_A$.
\end{pf}

Now, define the $A$-localized functor $\Pi ^A_0$ from $\deop \bcS $ to $\bcS $ by setting $\Pi ^A_0(\bcX )$ to be the sheaf associated to the presheaf
  $$
  U\mapsto [\Const (h_U),\bcX ]_A\; .
  $$
Then $\Pi ^A_0(\bcX )$ is canonically isomorphic to $\Pi _0(L_A(\bcX ))$, and the morphism $l$ induces the epimorphism $\Pi _0\to \Pi _0^A$, see Corollary 3.22 on page 94 in \cite{MorelVoevodsky}. As $L_A$ is monoidal,
  $$
  \begin{array}{rcl}
  \Pi ^A_0(\bcX \times \bcY )
  &=&
  \Pi _0(L_A(\bcX \times \bcY )) \\
  &=&
  \Pi _0(L_A(\bcX )\times L_A(\bcY )) \\
  &=&
  \Pi _0(L_A(\bcX ))\times \Pi _0(L_A(\bcY )) \\
  &=&
  \Pi ^A_0(\bcX )\times \Pi ^A_0(\bcY )\; .
  \end{array}
  $$
This gives that $\Pi _0^A$ preserves monoids and groups.

\begin{lemma}
\label{bloha}
For any monoid $\bcX $ in $\deop \bcS $, one has a canonical isomorphism
  $$
  \Pi _0(\bcX )^+\simeq \Pi _0(\bcX ^+)
  $$
in $\bcS $.
\end{lemma}

\begin{pf}
Since $\Gamma (\bcX ^+)=\Gamma (\bcX )^+$ and $\Pi _0(\bcX )^+$ are completions, one has the universal morphisms $\gamma $ from $\Gamma (\bcX ^+)$ to $\Pi _0(\bcX )^+$ and $\delta $ from $\Pi _0(\bcX )^+$ to $\Pi _0(\bcX ^+)$. Since $\Gamma (\bcX )=\bcX _0$, $\Gamma (\bcX ^+)=\bcX ^+_0$ and $\gamma \circ \Gamma (\nu _{\bcX })=\nu _{\Pi _0(\bcX )}\circ \Psi $, where $\nu $ stays for the corresponding canonical morphisms from the monoids to their completions, the two compositions $\bcX ^+_1\rightrightarrows \bcX ^+_0\stackrel{\gamma }{\to }\Pi _0(\bcX )^+$ coincide, which gives the universal morphism $\varepsilon $ from $\Pi _0(\bcX ^+)$ to $\Pi _0(\bcX )^+$. Since $\Psi $ is an epimorphism, and using the uniqueness of the appropriate universal morphisms, we show that $\delta $ and $\varepsilon $ are mutually inverse isomorphisms of groups in $\bcS $.
\end{pf}

Let $\CMon (\deop \bcS )$ be the category of commutative monoids in $\deop \bcS $. Suppose that all cofibrations in $\deop \bcS $ are symmetrizable, see \cite{GorchinskiyGuletskii}. Then the simplicial model structure on $\deop \bcS $ gives rise to a simplicial model structure on $\CMon (\deop \bcS )$, compatible with Bousfield localizations, see \cite{PavlovScholbach1}, \cite{PavlovScholbach2}, \cite{PavlovScholbach3}, \cite{White1} and \cite{White2}. A morphism in $\CMon (\deop \bcS )$ is a weak equivalence (respectively, fibration) if and only if it is a weak equivalence (respectively, fibration) in $\deop \bcS $, loc.cit. The classifying space functor $B$ is then a functor from the model category $\CMon (\deop \bcS )$ to the model category $\deop \bcS $. Lemma 2.35 on page 85 in \cite{MorelVoevodsky}, and the universality of a left localization of a model structure (see part (b) of the Definition 3.1.1 on pp 47 - 48 of \cite{Hirsch}), being applied to the functor $B$, yield a (simplicial) weak equivalence
  $$
  (B\circ L_A)(\bcX )\simeq (L_{S^1\wedge A}\circ B)(\bcX )\; ,
  $$
for any commutative monoid $\bcX $ in $\deop \bcS $.

\begin{lemma}
\label{Pi0&completion}
For any pointwise good commutative monoid $\bcX $ in $\deop \bcS $, one has a canonical isomorphism
  $$
  \Pi ^A_0(\bcX )^+\simeq
  \Pi ^A_0(\bcX ^+)\; .
  $$
in $\bcS $.
\end{lemma}

\begin{pf}
Since $\bcX $ is pointwise good, one has the isomorphism
  $$
  \bcX ^+\simeq (\Omega \circ \Ex \circ B)(\bcX )
  $$
in $\Ho $ by Lemma \ref{keylemma}, where $\Omega $ is the simplicial loop functor and $\Ex $ is the (pre-$A$-localized) fibrant replacement for simplicial sheaves. Applying $L_A$ we get the isomorphism
  $$
  L_A(\bcX ^+)\simeq
  L_A((\Omega \circ \Ex \circ B)(\bcX ))\; .
  $$
By Theorem 2.34 on page 84 in \cite{MorelVoevodsky},
  $$
  L_A((\Omega \circ \Ex \circ B)(\bcX ))\simeq
  (\Omega \circ \Ex \circ L_{S^1\wedge A})(B(\bcX ))\; .
  $$
Since $B\circ L_A\simeq L_{S^1\wedge A}\circ B$, we obtain the isomorphism
  $$
  L_A(\bcX ^+)\simeq (\Omega \circ \Ex \circ B)(L_A(\bcX ))
  $$
in $\Ho $. Let $\Phi =\Phi _{Mon}$ be the functor constructed in Lemma 1.1 on page 123 in \cite{MorelVoevodsky}, i.e. the cofibrant replacement functor in $\CMon (\deop \bcS )$. Since the morphism
  $$
  \Phi (L_A(\bcX ))\to L_A(\bcX )
  $$
is a weak equivalence in $\deop \bcS $, we get the isomorphism
  $$
  L_A(\bcX ^+)\simeq (\Omega \circ \Ex \circ B)(\Phi (L_A(\bcX )))
  $$
in $\Ho $. The monoid $\Phi (L_A(\bcX ))$ is termwise free. Therefore,
  $$
  (\Omega \circ \Ex \circ B)(\Phi (L_A(\bcX )))\simeq
  (\Phi (L_A(\bcX )))^+
  $$
by Lemma 1.2 on page 123 in \cite{MorelVoevodsky}. Applying $\Pi _0$ and using Lemma \ref{bloha}, we obtain the isomorphisms
  $$
  \begin{array}{rcl}
  \Pi ^A_0(\bcX ^+)
  &=&
  \Pi _0(L_A(\bcX ^+)) \\
  &=&
  \Pi _0((\Phi (L_A(\bcX )))^+) \\
  &=&
  \Pi _0(\Phi (L_A(\bcX )))^+ \\
  &=&
  \Pi _0(L_A(\bcX ))^+ \\
  &=&
  \Pi ^A_0(\bcX )^+ \\
  \end{array}
  $$
in the category of sheaves $\bcS $.
\end{pf}

\section{Chow monoids in Nisnevich sheaves}
\label{Chowsheaves}

Now we turn from homotopy algebra to algebraic geometry. Throughout all schemes will be separated by default. Let $k$ be a field, $\Sm $ the category of smooth schemes of finite type over $k$, and let $\goN $ be the category of all noetherian schemes over $k$, not necessarily of finite type. We are going to specialize the abstract material of the previous sections to the case when $\bcC $ is $\Sm $, the topology $\tau $ is the Nisnevich topology on $\bcC $, and $A$ is the affine line $\AF ^1$ over $k$.

The standard Yoneda construction gives the functor $h$ sending any scheme $X$ from $\goN $ to the functor $\Hom _{\goN }(-,X)$, and the same on morphisms. This is a functor to the category of sheaves in \'etale topology, and so in the Nisnevich one, see \cite{SGA4-2}, page 347, i.e. the Nisnevich topology is subcanonical. Composing $h$ with the constant functor $\Const =\Delta _?[0]$ from $\bcS $ to $\deop \bcS $ we obtain the embedding of $\goN $ into $\deop \bcS $. We identify the categories $\goN $ and $\SSets $ with their images under the corresponding embeddings into $\deop \bcS $.

The scheme $\Spec (k)$ is the terminal object in $\bcC $. The affine line $\AF ^1$ over $k$ is an interval in $\deop \bcS $ with two obvious morphisms $i_0$ and $i_1$ from $\Spec (k)$ to $\AF ^1$. As above, the interval $\AF ^1$ gives the natural cylinder and the corresponding notion of left homotopy on morphisms in $\deop \bcS $. The set of points on a simplicial sheaf $\bcX $ is the set $\Hom _{\deop \bcS }(\Spec (k),\bcX )$ of $k$-points on $\bcX $, and it coincides with the set $\bcX _0(k)$. The set of $\AF ^1$-path connected components on $k$-points is the set $\Hom _{\deop \bcS }(\Spec (k),\bcX )_{\AF ^1}$. If $\bcX $ is fibrant in $\bcM _{\AF ^1}$, then $\Hom _{\deop \bcS }(\Spec (k),\bcX )_{\AF ^1}$ can be identified with $[\Spec (k),\bcX ]_{\AF ^1}$.

Let $\bcX $ be a monoid in $\deop \bcS $. Its completion $\bcX ^+$ is a group object, so that $\Hom _{\deop \bcS }(\Spec (k),\bcX ^+)$ is a group in $\deop \bcS $. The morphism $\bcX \to \bcX ^+$ induces a map from $\Hom _{\deop \bcS }(\Spec (k),\bcX )$ to $\Hom _{\deop \bcS }(\Spec (k),\bcX ^+)$. By the universality of group completion, there exists a unique map from $\Hom _{\deop \bcS }(\Spec (k),\bcX )^+$ to $\Hom _{\deop \bcS }(\Spec (k),\bcX ^+)$ with the obvious commutativity.

\begin{lemma}
\label{muravei}
For a simplicial Nisnevich sheaf monoid $\bcX $, the canonical map from $\Hom _{\deop \bcS }(\Spec (k),\bcX )^+$ to $\Hom _{\deop \bcS }(\Spec (k),\bcX ^+)$ is bijective, and, repspectively, the map from $({\Hom _{\deop \bcS }(\Spec (k),\bcX )_{\AF ^1}})^+$ to $\Hom _{\deop \bcS }(\Spec (k),\bcX ^+)_{\AF ^1}$ is a surjection.
\end{lemma}

\begin{pf}
Since $\Spec (k)$ is Henselian, the set $\Hom _{\deop \bcS }(\Spec (k),\bcX ^+)$ is the quotient of the Cartesian square $\Hom _{\deop \bcS }(\Spec (k),\bcX )^2$. The set $\Hom _{\deop \bcS }(\Spec (k),\bcX )^+$ is also the quotient of the same Cartesian square. The maps from $\Hom $-sets to the sets of $\AF ^1$-path connected components are surjective.
\end{pf}

Next, let $K$ be a field extension of the ground field $k$, and let $\bcS _K$ be the category of set valued Nisnevich sheaves on the category $\bcC _K$ of smooth schemes over $K$. Let $\bcM _K$ be the injective model structure on the category $\deop \bcS _K$, obtained in the same way as the model structure $\bcM $ for the category $\deop \bcS $ over $\Spec (k)$. Let $f:\Spec (K)\to \Spec (k)$ be the morphism induced by the extension $k\subset K$, and let $f^*:\deop \bcS \to \deop \bcS _K$ be the scalar extension functor induced by sending schemes over $k$ to their fibred products with $\Spec (K)$ over $\Spec (k)$, and then using the fact that any sheaf is a colimit of representable ones. As the morphism $f$ is smooth, there are two standard adjunctions
  $$
  f_{\# } \dashv f^*\dashv f_*
  $$
for the functor $f^*$, see, for example, \cite{Morel2}.

\begin{lemma}
\label{ext}
For any field extension $K$ of the ground field $k$, one can choose the localization functors $L_{\AF ^1}$ and $L_{\AF ^1_K}$ in $\deop \bcS $ and $\deop \bcS _K$ respectively, to have a canonical isomorphism
  $$
  f^*L_{\AF ^1}\simeq L_{\AF ^1_K}f^*\; .
  $$

\end{lemma}

\begin{pf}
Let $\ExG\!\!_K$ be the fibrant replacement in $\deop \bcS _K$ obtained in the same way as $\ExG $ was constructed for $\deop \bcS $, see page 70 in \cite{MorelVoevodsky}. Let also $\Sing _K$ be the Suslin-Voevodsky endofunctor on $\deop \bcS _K$. Straightforward verifications show that $f^*\ExG \simeq \ExG\!\!_Kf^*$ and $f^*\circ \Sing \simeq \Sing _K\circ f^*$. Choose $L_{\AF ^1}$ (respectively, $L_{\AF ^1_K}$) to be the transfinite compositions of the functors $\ExG $ and $\Sing $ (respectively, $\ExG\!\!_K$ and $\Sing _K$) in $\deop \bcS $ (respectively, in $\deop \bcS _K$).
\end{pf}

We now need to refresh some things from \cite{SV-ChowSheaves}. For a scheme $X$ let $\bcC (X)$ be the free commutative monoid generated by points of $X$, and let $\bcZ (X)$ be the group completion of $\bcC (X)$. An algebraic cycle $\zeta $ is an element in $\bcZ (X)$. As such, $\zeta $ is a finite linear combination $\sum m_i\zeta _i$ of points $\zeta _i$ on $X$ with integral coefficients $m_i$. The cycle $\zeta $ is said to be effective if and only if $m_i\geq 0$ for all $i$. This is equivalent to say that $\zeta $ is an element of $\bcC (X)$.

The support $\supp (\zeta )$ of $\zeta $ is the union of the Zariski closures of the points $\zeta _i$ with the induced reduced structures on them. The correspondence between points on $X$ and the reduced irreducible closed subschemes of $X$ allows to consider algebraic cycles as linear combinations $Z=\sum m_iZ_i$, where $Z_i$ is the Zariski closure of the point $\zeta _i$, for each $i$. Then $\supp (Z)$ is the same thing as $\supp (\zeta )$. The points $\zeta _i$, or the corresponding reduced closed subschemes $Z_i$, are prime cycles on $X$. The dimension of a point in $X$ is the dimension of its Zariski closure in $X$. Let then $\bcC _r(X)$ be the submonoid in $\bcC (X)$ generated by points of dimension $r$, and, respectively, let $\bcZ _r(X)$ be the subgroup in the abelian group $\bcZ (X)$ generated by points of dimension $r$ in $X$.

Let $S$ be a Noetherian scheme, let $k$ be a field, and let
  $$
  P:\Spec (k)\to S
  $$
be a $k$-point of $S$. Recall that a fat point of $S$ over $P$ is an ordered pair $(P_0,P_1)$ of two morphisms of schemes
  $$
  P_0:\Spec (k)\to \Spec (R)\quad \hbox{and}\quad P_1:\Spec (R)\to S\; ,
  $$
where $R$ is a discrete valuation ring with the residue field $k$, such that
  $$
  P_1\circ P_0=P\; ,
  $$
the image of $P_0$ is the closed point of $\Spec (R)$, and $P_1$ sends the generic point $\Spec (R_{(0)})$ into the generic point of $S$.

Let $X\to S$ be a scheme of finite type over $S$, and let
  $$
  Z\to X
  $$
be a closed subscheme in $X$. Let $R$ be a discrete valuation ring, $D=\Spec (R)$, and let
  $$
  f:D\to S
  $$
be an arbitrary morphism of schemes from $D$ to $S$. Let also
  $$
  \eta =\Spec (R_{(0)})
  $$
be the generic point of $D$,
  $$
  X_D=X\times _SD\; ,\quad Z_D=Z\times _SD\quad \hbox{and}\quad
  Z_{\eta }=Z\times _S\eta \; .
  $$
Then there exists a unique closed embedding
  $$
  Z'_D\to Z_D\; ,
  $$
such that its pull-back $Z'_{\eta }\to Z_{\eta }$, with respect to the morphism $Z_{\eta }\to Z_D$, is an isomorphism, and the composition
  $$
  Z'_D\to Z_D\to D
  $$
is a flat morphism of schemes, see Proposition 2.8.5 in \cite{EGAIV(2)}.

In particular, we have such a ``platification" if $(P_0,P_1)$ is a fat point over the $k$-point $P$ and $f=P_1$. Let then $X_P$ be the fibre of the morphism $X_D\to D$ over the point $P_0$,
  $$
  Z_P=Z_D\times _{X_D}X_P
  \quad \hbox{and}\quad
  Z'_P=Z'_D\times _{Z_D}Z_P\; .
  $$
Since the closed subscheme $Z'_D$ of $X_D$ is flat over $D$, one can define the pull-back
  $$
  (P_0,P_1)^*(Z)
  $$
of the closed subscheme $Z$ to the fibre $X_P$ of the morphism $X\to S$, with regard to the fat point $P$, as the cycle associated to the closed embedding $Z'_P\to X_P$ in the standard way (consult \cite[1.5]{Fulton} for what ``the standard way" means).

In particular, if $Z$ is a prime cycle on $X$, then we have the pull-back cycle $(P_0,P_1)^*(\zeta )$ on $X_P$. Extending by linearity we obtain a pull-back homomorphism
  $$
  (P_0,P_1)^*:\bcZ (X)\to \bcZ (X_P)\; .
  $$
Following \cite{SV-ChowSheaves}, we say that an algebraic cycle $\zeta =\sum m_i\zeta _i$ on $X$ is a relative cycle on $X$ over $S$ if the images of the points $\zeta _i$ under the morphism $X\to S$ are the generic points of the scheme $S$, and, for any $k$-point $P$ on $S$, and for any two fat points extending $P$, the pull-backs of the cycle $Z=\bar \zeta $ to $X_P$, with regard to these two fat points, coincide, see Definition 3.1.3 in loc.cit.

Notice that any cycle, which is flat over $S$, is a relative cycle for free. But not any relative cycle on $X/S$ is flat. This is why we need the ``platification" above.

Let $\bcC _r(X/S)$ be the abelian submonoid in $\bcC (X)$ generated by relative cycles of relative dimension $r$ over $S$. It is important that whenever $S$ is a regular Noetherian scheme and $X$ is of finite type over $S$, then $\bcC _r(X/S)$ is a free commutative monoid generated by closed integral subschemes in $X$ which are equidimensional of dimension $r$ over $S$, see Corollary 3.4.6 on page 40 in \cite{SV-ChowSheaves}. Let also
  $$
  \bcZ _r(X/S)=\bcC _r(X/S)^+
  $$
be the group completion of the monoid $\bcC _r(X/S)$.

Now, fix a Noetherian reduced scheme $T$, and let $\goN $ be the category of Noetherian schemes over $T$. Let $X\to T$ be a scheme of finite type over $T$. For any object $S\to T$ in $\goN $ let
  $$
  \bcZ _r(X/T)(S)=\bcZ _r(X\times _TS/S)\; .
  $$
If $f:S'\to S$ is a flat morphism of Noetherian schemes over $T$, the induced morphism $\id _X\times _Tf:X\times _TS'\to X\times _TS$ is also flat, and one has the standard flat pull-back homomorphism
  $$
  (\id _X\times _Tf)^*:\bcZ _r(X\times _TS/S)\to
  \bcZ _r(X\times _TS'/S')\; .
  $$
If $f$ is not flat, then the situation is more difficult.

However, if $T$ is a regular scheme, due to the Suslin-Voevodsky's definition of relative cycles given above, the correct pull-back exists for any morphism $f$, see Proposition 3.3.15 in \cite{SV-ChowSheaves}.

This all aggregates, when $T$ is a regular scheme, into the presheaf $\bcZ _r(X/T)$, which is nothing else but the sectionwise completion of the presheaf $\bcC _r(X/T)$ of relative effective cycles on the category $\goN $.

Since now we will assume that $T$ is regular of characteristic $0$. Let $X\to T$ be a projective scheme over $T$, and fix a closed embedding
  $$
  X\to \PR ^m_T
  $$
over $T$. If $Z$ is a relative equidimensional cycle on $X\times _TS/S$, its pullback $Z_P$ to the fibre $X_P$ of the morphism $X\times _TS\to S$ over a point $P$ on $S$ has its degree $\deg (Z_P)$, computed with regard to the induce embedding of $X\times _TS$ into $\PR ^m_S$ over $S$. Since $Z$ is a relative cycle, the degree $\deg (Z_P)$ is locally constant on $S$, see Proposition 4.4.8 in \cite{SV-ChowSheaves}. It follows that, if $S$ is connected, then $\deg (Z_P)$ does not depend on $P$, see Corollary 4.4.9 in loc.cit. Therefore, the degree of $Z$ over $S$ is correctly defined, and we may consider a subpresheaf
  $$
  \bcC _{r,d}(X/T)\subset \bcC _r(X/T)\; ,
  $$
whose sections on $S$ are relative cycles of degree $d$ on $X\times _TS/S$.

The integer $d$ is non-negative, and there is only one cycle in the set $\bcC _{r,0}(X/T)(S)$, namely the cycle $0$ whose coefficients are all zeros. The grading by degrees gives the obvious structure of a graded monoid on the presheaf $\bcC _r(X/T)$ whose neutral element is the cycle $0$ sitting in $\bcC _{r,0}(X/T)(S)$.

It follows from the results in \cite{SV-ChowSheaves} (see also \cite{Kollar}) that the presheaves $\bcC _{r,d}(X/T)$ are representable by a scheme
  $$
  C_{r,d}(X/T)\; ,
  $$
the so-called {\it Chow scheme} of effective relative cycles of relative dimension $r$ and degree $d$ over $T$. This Chow scheme is projective over $T$, i.e. there exist a structural morphism from $C_{r,d}(X/T)$ to $T$, and a closed embedding of $C_{r,d}(X/T)$ into $\PR ^N_T$ over $T$, arising from the representability above.

Notice that the Chow sheaves are representable because $T$ is a regular scheme of characteristic zero. If $T$ would be of positive characteristic, then only $h$-representability takes place, see \cite{SV-ChowSheaves}. If $T$ is the spectrum of an algebraically closed field of characteristic zero, then $C_{r,d}(X/T)$ is the classical Chow scheme of effective $r$-cycles of degree $d$ on $X$.

Since now we will assume that $T=\Spec (k)$, where $k$ is a field of characteristic zero, and systematically drop the symbol $/T$ from the notation. According to our convention to identify schemes and the corresponding representable sheaves, we will write $C_{r,d}(X)$ instead of $h_{C_{r,d}(X)}$. Certainly, the latter sheaf $C_{r,d}(X)$ is isomorphic to, and should be identified with the sheaf $\bcC _{r,d}(X)$.

Let
  $$
  C_r(X)=\coprod _{d\geq 0}C_{r,d}(X)\; ,
  $$
where the coproduct is taken in the category $\bcS $, not in $\goN $. Such defined $C_r(X)$ is also a coproduct in $\bcP $. If we would consider the coproduct of all Chow schemes $C_{r,d}(X)$ in $\goN $ first, and then embed it into $\bcS $ by the Yoneda embedding, that would be a priori a different sheaf, as Yoneda embedding in general does not commute with coproducts. However, the canonical morphism from the above sheafification to this second sheaf is an isomorphism on the Henselizations of the local rings at points of varieties over $k$. Therefore, the two constructions are actually isomorphic in $\bcS $. This also gives that the coproduct of $C_{r,d}(X)$, for all $d\geq 0$, in $\goN $ represents $\bcC _r(X)$.

Identifying $\bcS $ with its image in $\deop \bcS $ under the functor $\Const $, we consider $C_r(X)$ as the graded {\it Chow monoid} in the category of simplicial sheaves on the smooth Nisnevich site over $\Spec (k)$. The completion $C_r(X)^+$ of $C_r(X)$ in $\bcS $ is the sheafification of the completion of $C_r(X)$ as a presheaf. The latter is sectionwise.

Let $\bcO ^h_{P,Y}$ be the Henselization of the local ring $\bcO _{P,Y}$ of a smooth algebraic scheme $Y$ over $k$ at a point $P\in Y$. Since $\bcO _{P,Y}$ is a regular Noetherian ring, the ring $\bcO _{P,Y}^h$ is regular and Noetherian too. As we mentioned above, the set
  $$
  \bcC _r(X\times _{\Spec (k)}\Spec (\bcO ^h_{P,Y})/\Spec (\bcO ^h_{P,Y}))
  $$
is a free commutative monoid generated by closed integral subschemes in the scheme $X\times _{\Spec (k)}\Spec (\bcO _{P,Y}^h)$, which are equidimensional of dimension $r$ over $\Spec (\bcO _{P,Y}^h)$, by Corollary 3.4.6 in \cite{SV-ChowSheaves}. Then we see that the monoid $C_r(X)$ is pointwise free, and hence it is a pointwise cancellation monoid in the category $\bcS $. It follows also that the Chow monoid $C_r(X)$ is pointwise good in the category $\deop \bcS $, and the canonical morphism from $C_r(X)$ to $C_r(X)^+$ is a monomorphism.

Let $K$ be a finitely generated field extension of $k$. Since $\Spec (K)$ is Henselian, $C_r(X)^+(K)$ is the same as the group completion $(C_r(X)(K))^+$. On the other hand, the same group $C_r(X)^+(K)$ can be also identified with the group of morphisms from $\Spec (K)$ to $C_r(X)^+$, in the category of simplicial Nisnevich sheaves $\deop \bcS $.

Let $d_0$ be the minimal degree of positive $r$-cycles on $X$, where the degree is computed with regard to the fixed embedding of $X$ into $\PR ^m$. Choose and fix a positive $r$-cycle $Z_0$ with $\deg (Z_0)=d_0$. For any natural number $d$ the $d$-multiple $dZ_0$ is an effective dimension $r$ degree $dd_0$ cycle on $X$. This gives a morphism $\alpha $ from $\NN $ to $C_r(X)$ sending $1$ to $Z_0$. Since $C_r(X)$ is the coproduct of $C_{r,dd_0}(X)$, for all $d\geq 0$, we also have the obvious morphism $f$ from $C_r(X)$ to $\NN $, such that $f\circ \alpha =\id _{\NN }$. In other words, $Z_0$ gives the structure of a pointed graded monoid on $C_r(X)$. Automatically, we obtain the connective Chow monoid $C_r^{\infty }(X)$ associated to $C_r(X)$. By Lemma \ref{goraimysh}, we also have the canonical isomorphism of group objects
  $$
  C_r(X)^+\simeq C_r^{\infty }(X)^+\times \ZZ
  $$
in $\deop \bcS $. The sheaf $C_r^{\infty }(X)$ can be also understood as the ind-scheme arising from the chain of Chow schemes
  $$
  C_{r,0}(X)\subset C_{r,d_0}(X)\subset C_{r,2d_0}(X)\subset \ldots \subset C_{r,dd_0}(X)
  \subset \ldots
  $$
induced by the cycle $Z_0$ of degree $d_0$ on $X$.

As the category $\bcC $ is Noetherian, the category $\bcS $ is exhaustive. Since $C_r(X)$ is a pointwise cancellation monoid in $\bcS $, and the latter category is exhaustive, $C_r^{\infty }(X)$ is a pointwise cancellation monoid in $\bcS $ too. Then both monoids, $C_r(X)$ and $C_r^{\infty }(X)$ are pointwise good monoids in the category $\deop \bcS $. Moreover, the canonical morphism from $C_r^{\infty }(X)$ to $C_r^{\infty }(X)^+$ is a monomorphism in $\bcS $ and in $\deop \bcS $.

\section{Rational equivalence as $\AF ^1$-path connectivity}
\label{rateq-ratcon}

For any algebraic scheme $X$ over $k$ let $CH_r(X)$ be the Chow group of $r$-dimensional algebraic cycles modulo rational equivalence on $X$. In this section we prove our main theorem and deduce three corollaries, which give something close to the desired effective interpretation of Chow groups in terms of $\AF ^1$-path connectivity on loop spaces of classifying spaces of the Chow monoid $C_r^{\infty }(X)$. We leave it for the reader to decide which of the obtained three isomorphisms is more useful for understanding of Chow groups.

\begin{theorem}
\label{main1}
Let $X$ be a projective algebraic variety with a fixed embedding into a projective space over $k$. For any finitely generated field extension $K$ of the ground field $k$, there is a canonical isomorphism
  $$
  CH_r(X_K)\simeq
  \Pi _0^{\AF ^1}(C_r(X)^+)(K)\; .
  $$
\end{theorem}

\medskip

\begin{pf}
Without loss of generality, we may assume that $d_0=1$. Consider the obvious commutative diagram
  $$
  \diagram
  \Hom _{\deop \bcS }
  (\Spec (k),C_r(X)^+)_{\AF ^1}
  \ar[r]^-{(l^+)_*} &
  \Hom _{\deop \bcS }
  (\Spec (k),L_{\AF ^1}C_r(X)^+)_{\AF ^1} \\ \\
  {\Hom _{\deop \bcS }
  (\Spec (k),C_r(X))_{\AF ^1}}^+
  \ar[uu]_-{\alpha } \ar[r]^-{(l_*)^+} &
  {\Hom _{\deop \bcS }
  (\Spec (k),L_{\AF ^1}C_r(X))_{\AF ^1}}^+
  \ar[uu]_-{\beta }
  \enddiagram
  $$
where $l=l_{C_r(X)}$. Since $L_{\AF ^1}C_r(X)^+$ is $\AF ^1$-local, the group in the top right corner is canonically isomorphic to the group $\Pi _0^{\AF ^1}(C_r(X)^+)(k)$. By the same reason, the group in the bottom right corner is canonically isomorphic to the group $\Pi _0^{\AF ^1}(C_r(X))(k)^+$. Since $\Spec (k)$ is Henselian, the latter group is nothing but the group $\Pi _0^{\AF ^1}(C_r(X))^+(k)$. Then Lemma \ref{Pi0&completion} gives that $\beta $ is an isomorphism.

Let $q_0:\Spec (k)\to C_r(X)$ and $q_1:\Spec (k)\to C_r(X)$
be two $k$-points on $C_r(X)$, and suppose $q_0$ is connected to $q_1$ by an $\AF ^1$-path $H:\AF ^1\to L_{\AF ^1}C_r(X)$ on $L_{\AF ^1}C_r(X)$. For any $d$ let $(C_r(X))_d$ be the coproduct $\coprod _{i=0}^dC_{r,i}(X)$. Then $(C_r(X))_d$ is canonically embedded into the coproduct $(C_r(X))_{d+1}$. Consider the chain of the embeddings
  $$
  (C_r(X))_0\subset \ldots \subset (C_r(X))_d\subset (C_r(X))_{d+1}
  \subset \ldots
  $$
Applying Proposition 4.4.4 on page 77 in \cite{Hirsch} (see also Remark \ref{filin} in Section \ref{localization}) we see that $L_{\AF ^1}C_r(X)$ is canonically isomorphic to the colimit of the chain of the embeddings
  $$
  L_{\AF ^1}((C_r(X))_0)\subset \ldots \subset
  L_{\AF ^1}((C_r(X))_d)\subset L_{\AF ^1}((C_r(X))_{d+1})
  \subset \ldots
  $$
Since $\AF ^1$ is a compact object in the category $\deop \bcS $, it follows that the homotopy $H$ factorizes through $L_{\AF ^1}((C_r(X))_d)$, for some degree $d$. If $Z_0$ is a degree $1$ algebraic cycle of dimension $r$ on $X$, then $Z_0$ induces the corresponding embeddings
  $$
  C_{r,0}(X)\subset \ldots \subset C_{r,d}(C)\; .
  $$
This gives the epimorphism from the coproduct $(C_r(X))_d$ onto $C_{r,d}(X)$. Composing the homotopy $H:\AF ^1\to L_{\AF ^1}((C_r(X))_d)$ with the induced morphism from $L_{\AF ^1}((C_r(X))_d)$ to $L_{\AF ^1}(C_{r,d}(X))$, we obtain the homotopy
  $$
  H:\AF ^1\to L_{\AF ^1}(C_{r,d}(X))\; .
  $$
Since $X$ is proper and of finite type over $k$, Proposition 6.2.6 in \cite{AsokMorel} gives that the points $q_0$ and $q_1$ are $\AF ^1$-chain connected, and so $\AF ^1$-path connected on $C_{r,d}(X)$. It means that the map
  $$
  l_*:\Hom _{\deop \bcS }(\Spec (k),C_r(X))_{\AF ^1}
  \to
  \Hom _{\deop \bcS }
  (\Spec (k),L_{\AF ^1}C_r(X))_{\AF ^1}
  $$
is injective. Since $l_*$ is surjective by Lemma \ref{triton}, it is bijective. Then $(l^+)_*$ is an isomorphism as well.

Since $\beta $ and $(l^+)_*$ are isomorphisms, and $\alpha $ is an epimorphism by Lemma \ref{muravei}, we see that $\alpha $ is an isomorphism, and then all the maps in the commutative square above are isomorphisms.

Let now $A$ and $A'$ be two $r$-dimensional algebraic cycles on $X$. If $A$ is rationally equivalent to $A'$ on $X$, there exists an effective relative cycle $Z$ on the scheme $X\times _{\Spec (k)}\AF ^1/\AF ^1$ of relative dimension $r$, and an effective dimension $r$ algebraic cycle $B$ on $X$, such that
  $$
  Z(0)=A+B\quad \hbox{and}\quad Z(1)=A'+B
  $$
on $X$. Let $h_Z$ and $h_{B\times \AF ^1}$ be two regular morphisms from $\AF ^1$ to the Chow scheme $C_r(X)$ over $\Spec (k)$ corresponding to the relative cycles $Z$ and $B\times \AF ^1$ on $X\times _{\Spec (k)}\AF ^1/\AF ^1$ respectively. Let
  $$
  h:\AF ^1\to C_r(X)\times C_r(X)
  $$
be the product of $h_Z$ and $h_{B\times \AF ^1}$ in the category $\deop \bcS $. Let
  $$
  H:\AF ^1\to C_r(X)\times C_r(X)\to C_r(X)^+\; ,
  $$
be the composition of $h$ and the morphism from $C_r(X)\times C_r(X)$ to the completion $C_r(X)^+$, in $\deop \bcS $. Then $H_0=A$ and $H_1=A'$, where $H_0$ and $H_1$ are the precompositions of $H$ with $i_0$ and $i_1$ respectively. It means that the cycles $A$ and $A'$ are $\AF ^1$-path connected on $C_r(X)^+$.

Vice versa, suppose we have a morphism
  $$
  H:\AF ^1\to C_r(X)^+
  $$
in $\bcS $, and let $H_0$ and $H_1$ be the compositions of $H$ with $i_0$ and $i_1$ respectively. Since $\Spec (k)$ is Henselian, $H_0$ is represented by two morphisms $H_{0,1}$ and $H_{0,2}$ from $\Spec (k)$ to $C_r(X)$. Similarly, $H_1$ is represented by two morphisms $H_{1,1}$ and $H_{1,2}$ from $\Spec (k)$ to $C_r(X)$. Since $\alpha $ is an isomorphism and $H_0$ is $\AF ^1$-path connected to $H_1$, it follows that there exist two morphisms $f$ and $G$ from $\Spec (k)$ to $C_r(X)$, such that $H_{0,1}+F$ is $\AF ^1$-path connected to $H_{0,2}+G$ and $H_{1,1}+F$ is $\AF ^1$-path connected to $H_{1,2}+G$ on $C_r(X)$. In terms of algebraic cycles on $X$, it means that the effective $r$-cycle $H_{0,1}+F$ is rationally equivalent to the effective $r$-cycle $H_{0,2}+G$, and, similarly, the cycle $H_{1,1}+F$ is rationally equivalent to $H_{1,2}+G$ on $C_r(X)$. Then the cycle $H_0=H_{0,1}-H_{0,2}$ is rationally equivalent to the cycle $H_1=H_{1,1}-H_{1,2}$ on $C_r(X)$.

Thus, the Chow group $CH_r(X)$ is isomorphic to $\Hom _{\deop \bcS }(\Spec (k),C_r(X)^+)_{\AF ^1}$, i.e. the group in the top left corner of the diagram above. Since, moreover, $(l^+)_*$ is an isomorphism, and the group in the top right corner is canonically isomorphic to $\Pi _0^{\AF ^1}(C_r(X)^+)(k)$, we obtain the required isomorphism in case when $L$ is the ground field $k$.

To prove the theorem for an arbitrary finitely generated field extension $K$ of the ground field $k$, we observe that $f^*C_r(X)$ is $C_r(X_K)$, whence
  $$
  f^*L_{\AF ^1}C_r(X)^+=
  L_{\AF ^1_K}f^*C_r(X)^+=
  L_{\AF ^1_K}C_r(X_K)^+
  $$
by Lemma \ref{ext}. Therefore,
  $$
  \begin{array}{rcl}
  \Pi _0^{\AF ^1}(C_r(X)^+)(K)
  &=&
  \Pi _0(L_{\AF ^1}C_r(X)^+)(K) \\
  &=&
  f^*\Pi _0(L_{\AF ^1}C_r(X)^+)(K) \\
  &=&
  \Pi _0(f^*L_{\AF ^1}C_r(X)^+)(K)\\
  &=&
  \Pi _0(L_{\AF ^1_K}C_r(X_K)^+)(K)\\
  &=&
  \Pi _0^{\AF ^1_K}(C_r(X_K)^+)(K) \\
  &\simeq &
  CH_r(X_K)_0\; .
  \end{array}
  $$
\end{pf}

\begin{remark}
{\rm Lemma \ref{Pi0&completion} provides that the monoidal completion in Theorem \ref{main1} can be taken before or after computing the $\AF ^1$-connected component functor. Since the monoidal completion is sectionwise on stalks, we obtain the canonical isomorphisms
  $$
  \begin{array}{rcl}
  CH_r(X)
  &\simeq &
  \Pi _0^{\AF ^1}(C_r(X)^+)(k) \\
  &\simeq &
  \Pi _0^{\AF ^1}(C_r(X))^+(k) \\
  &\simeq &
  \Pi _0^{\AF ^1}(C_r(X))(k)^+\; .
  \end{array}
  $$
}
\end{remark}

The embedding $X\hra \PR ^m$ gives the degree homomorphism from $CH_r(X)$ to $\ZZ $. Let $CH_r(X)_0$ be its kernel, i.e. the Chow group of degree $0$ cycles of dimension $r$ modulo rational equivalence on $X$. Then,
  $$
  CH_r(X)\simeq CH_r(X)_0\times \ZZ \; .
  $$
Let $Z_0$ be a positive $r$-cycle of minimal degree on $X$. As we have seen above, this gives the structure of a pointed graded cancellation monoid on $C_r(X)$, and $C_r^{\infty }(X)$ is a cancelation monoid too.

\begin{corollary}
\label{main2}
In terms above,
  $$
  CH_r(X_K)_0\simeq
  \Pi _0^{\AF ^1}(C_r^{\infty }(X)^+)(K)\; .
  $$
\end{corollary}

\begin{pf}
By Lemma \ref{goraimysh},
  $$
  C_r(X)^+\simeq C_r^{\infty }(X)^+\times \ZZ
  $$
Since the functor $\Pi _0^{\AF ^1}$ is monoidal and $\Pi _0^{\AF ^1}(\ZZ )=\ZZ $, we get the formula
  $$
  \Pi _0^{\AF ^1}(C_r(X)^+)\simeq
  \Pi _0^{\AF ^1}(C_r^{\infty }(X)^+)\times \ZZ \; .
  $$
Then apply Theorem \ref{main1} and the isomorphism $CH_r(X)\simeq CH_r(X)_0\times \ZZ $.
\end{pf}

\begin{warning}
{\rm If $CH_r(X)_0\simeq \Pi _0^{\AF ^1}(C_r^{\infty }(X))(k)^+=0$, it does not imply that the monoid $\Pi _0^{\AF ^1}(C_r^{\infty }(X))(k)$ vanishes, as this monoid is by no means a pointwise cancellation monoid. One of the reasons for that is that the Chow schemes $C_{r,d}(X)$ can have many components over $k$.
}
\end{warning}

\begin{corollary}
\label{main2.5}
In terms above,
  $$
  CH_r(X_K)_0\simeq
  \Pi _0^{\AF ^1}
  (\Omega \Ex BC_r^{\infty }(X))(K)\; .
  $$
\end{corollary}

\begin{pf}
The Chow monoid $C_r^{\infty }(X)$ is pointwise good. Lemma \ref{keylemma} gives an isomorphism
  $$
  C_r^{\infty }(X)^+\simeq \Omega \Ex BC_r^{\infty }(X)
  $$
in $\Ho $, whence
  $$
  \Pi _0^{\AF ^1}(C_r^{\infty }(X)^+)\simeq
  \Pi _0^{\AF ^1}(\Omega \Ex BC_r^{\infty }(X))\; .
  $$
Corollary \ref{main2} completes the proof.
\end{pf}

Recall that, for a pointed simplicial Nisnevich sheaf $(\bcX ,x)$, its motivic, i.e. $\AF ^1$-fundamental group $\Pi _1^{\AF ^1}(\bcX ,x)$ is, by definition, the Nisnevich sheaf associated to the presheaf sending a smooth scheme $U$ to the set $[S^1\wedge U_+,(\bcX ,x)]_{\AF ^1}$, where the symbol $[-,-]_{\AF ^1}$ stays now for the sets of morphisms in the pointed homotopy category $\Ho _{\AF ^1}$, see \cite{MorelVoevodsky} or \cite{AsokMorel}. Similarly, one can define, for a pointed simplicial Nisnevich sheaf $(\bcX ,x)$, the fundamental group $\Pi _1^{S^1\wedge \AF ^1}(\bcX ,x)$, where $\AF ^1$ is pointed at any $k$-rational point on it. This is the Nisnevich sheaf associated to the presheaf sending a smooth scheme $U$ to the set
  $$
  [S^1\wedge U_+,(\bcX ,x)]_{S^1\wedge \AF ^1}\; ,
  $$
where the symbol $[-,-]_{S^1\wedge \AF ^1}$ stays for the sets of morphisms in the pointed homotopy category $\Ho _{S^1\wedge \AF ^1}$.

\begin{lemma}
\label{onemorelemma}
Let $\bcX $ be a pointwise good simplicial sheaf monoid. Then, for a scheme $U$,
  $$
  \Pi _0^{\AF ^1}(\bcX ^+)(U)\simeq
  \Pi _1^{S^1\wedge \AF ^1}(B\bcX )(U)\; .
  $$
\end{lemma}

\begin{pf}
Since $\bcX $ is pointwise good, there is a isomorphism between $\bcX ^+$ and $\Omega \Ex B\bcX $ in the homotopy category $\Ho $, by Lemma \ref{keylemma}. Since the classifying space $B\bcX $ is pointed connected, the canonical morphism
  $$
  L_{\AF ^1}\Omega \Ex B\bcX \to
  \Omega \Ex L_{S^1\wedge \AF ^1}B\bcX
  $$
is a simplicial (i.e. pre-$\AF ^1$-localized) weak equivalence and
  $$
  \Omega \Ex L_{S^1\wedge \AF ^1}B\bcX
  $$
is $\AF ^1$-local by Theorem 2.34 on page 84 in loc.cit. This allows us to make the following identifications:
  $$
  \begin{array}{rcl}
  \Pi _0^{\AF ^1}(\bcX ^+)(U)
  &\simeq &
  [U,\bcX ^+]_{\AF ^1} \\
  &\simeq &
  [U,\Omega \Ex B\bcX ]_{\AF ^1} \\
  &\simeq &
  [U,L_{\AF ^1}\Omega \Ex B\bcX ]_{\AF ^1} \\
  &\simeq &
  [U,\Omega \Ex L_{S^1\wedge \AF ^1}B\bcX ]_{\AF ^1} \\
  &\simeq &
  [U,\Omega \Ex L_{S^1\wedge \AF ^1}B\bcX ] \\
  &\simeq &
  [S^1\wedge U_+,L_{S^1\wedge \AF ^1}B\bcX ] \\
  &\simeq &
  [S^1\wedge U_+,B\bcX ]_{S^1\wedge \AF ^1} \\
  &\simeq &
  \Pi _1^{S^1\wedge \AF ^1}(B\bcX )(U)\; .
  \end{array}
  $$

\end{pf}

\begin{corollary}
\label{main3}
In terms above,
  $$
  CH_r(X_K)_0\simeq
  \Pi _1^{S^1\wedge \AF ^1}(BC_r^{\infty }(X))(K)\; .
  $$
\end{corollary}

\begin{pf}
This is a straightforward consequence of Corollary \ref{main2} and Lemma \ref{onemorelemma}.
\end{pf}

\begin{example}
{\rm Let $X$ be a nonsingular projective surface over $k$, where $k$ is algebraically closed of characteristic zero. Assume that $X$ is of general type and has no transcendental second cohomology group, i.e. the cycle class map from $CH^1(X)$ to the second Weil cohomology group $H^2(X)$ is surjective. In that case the irregularity of $X$ is zero. Bloch's conjecture predicts that $CH_0(X)=\ZZ $. In other words, any two closed points on $X$ are rationally equivalent to each other. Fixing a point on $X$ gives the Chow monoid $C_0^{\infty }(X)$, which is nothing else but the the infinite symmetric power $\Sym ^{\infty }(X)$ of the smooth projective surface $X$. By Corollary \ref{main2.5}, Bloch's conjecture holds for $X$ if and only if all $k$-points on the motivic space $L_{\AF ^1}\Omega \Ex B\Sym ^{\infty }(X)$ are $\AF ^1$-path connected. Bloch's conjecture holds, for example, for the classical Godeaux surfaces, \cite{Sur les zero-cycles}, and for the Catanese and Barlow surfaces, see \cite{Barlow2} and \cite{VoisinCataneseBarlowSurfaces}.
}
\end{example}


The above vision of Chow groups should be compared with the results of Friedlander, Lawson, Lima-Filho, Mazur and others, who considered topological (i.e. not motivic) homotopy completions of Chow monoids working over $\CC $, see \cite{FriedlanderMazur} and \cite{Lima-Filho}. A nice survey of this topic, containing many useful references, is the article \cite{Lawson}. The topological homotopy completions of Chow monoids are helpful to understand algebraic cycles modulo algebraic equivalence relation, i.e. the groups $A_r(X)$ of algebraically trivial $r$-cycles cannot be catched by the topological methods. In contrast, the motivic, i.e. $\AF ^1$-homotopy completions of Chow monoids, considered above, can give the description of $A_r(X)$, working over an arbitrary ground field of characteristic zero, as the previous examples show. Theorem \ref{main1} also suggests that the {\it motivic Lawson homology} groups can be defined by the formula
  $$
  L_rH^{\bcM }_n(X)=\Pi _{n-2r}^{\AF ^1}(C_r(X)^+)(k)\; .
  $$


\medskip

\section{Appendix: homotopical algebra}

For the convenience of the reader, we collect here the needed extractions from homotopical algebra. Let first $\bcC $ be a symmetric monoidal category with product $\otimes $ and unit $\uno $. The monoidal product $\otimes $ is called to be closed, and the category $\bcC $ is called closed symmetric monoidal, if the product $\otimes :\bcC \times \bcC \to \bcC $ is so-called adjunction of two variables, i.e. there is bifunctor $\cHom $ and two functorial in $X$, $Y$, $Z$ bijections
  $$
  \Hom _{\bcC }(X,\cHom (Y,Z))\simeq
  \Hom _{\bcC }(X\otimes Y,Z)\simeq
  \Hom _{\bcC }(Y,\cHom (X,Z))\; .
  $$
If $\bcC $ has a model structure $\bcM $ in it, an adjunction of two variables on $\bcC $ is called Quillen adjunction of two variables, or Quillen bifunctor, if, for any two cofibrations $f:X\to Y$ and $f':X'\to Y'$ in $\bcM $ the push-out product
  $$
  f\square f':(X\otimes Y')\coprod _{X\otimes X'}
  (Y\otimes X')\to Y\wedge Y'
  $$
is also a cofibration in $\bcM $, which is trivial if either $f$ and $f'$ is. The model category $(\bcC ,\bcM )$ is called closed symmetric monoidal model category if $\otimes $ is a Quillen bifunctor and the following extra axiom holds. If $q:Q\uno \to \uno $ is a cofibrant replacement for the unit object $\uno $, then the morphisms $q\wedge \id :Q\uno \wedge X\to \uno \wedge X$ and $\id \wedge q:X\wedge Q\uno \to X\wedge \uno $ are weak equivalences for all cofibrant objects $X$. If we consider the cartesian product $\bcM \times \bcM $ of the model structure $\bcM $ as a model structure on the cartesian product $\bcC \times \bcC $, then $\otimes $ and $\cHom $ induce left derived functor $\otimes ^L$ from $Ho(\bcC \times \bcC )$ to $Ho(\bcC )$, and right derived functor $R\cHom $ from $Ho(\bcC \times \bcC )$ to $Ho(\bcC )$. It is well known that passing to localization commutes with products of categories,
so that we have the equivalence between $Ho(\bcC \times \bcC )$ and $Ho(\bcC )\times Ho(\bcC )$. This gives the left derived functor
  $$
  \otimes ^L:Ho(\bcC )\times Ho(\bcC )\to Ho(\bcC )
  $$
and the right derived functor
  $$
  R\cHom :Ho(\bcC )\times Ho(\bcC )\to Ho(\bcC )\; .
  $$
As it was shown in \cite{Hovey}, the left derived $\otimes ^L$ and the right derived $R\cHom $ give the structure of a closed symmetric monoidal category on the homotopy category $Ho (\bcC )$. Since we assume that all objects in $\bcC $ are cofibrant in $\bcM $, it is easy to see that the canonical functor from $\bcC $ to $Ho (\bcC )$ is monoidal.

An important particular case is when the symmetric monoidal product $\otimes $ is given by the categorical product in $\bcC $, i.e. when $\bcC $ is the cartesian symmetric monoidal category. Since $Ho(\bcC )$ admits products, and products in $\bcC $ are preserved in $Ho (\bcC )$,
for any three objects $X$, $Y$ and $Z$ in $\bcC $ one has the canonical isomorphism
  $$
  [X,Y]\times [X,Z]\simeq [X,Y\times Z]\; .
  $$

Let now $\bcC $ be a left proper cellular simplicial model category with model structure $\bcM =(W,C,F)$ in it, let $I$ and $J$ be the sets of, respectively, generating cofibrations and generating trivial cofibrations in $\bcC $, and let $S$ be a set of morphisms in $\bcC $. For simplicity we will also be assuming that all objects in $\bcC $ are cofibrant, which will always be the case in applications. An object $Z$ in $\bcC $ is called $S$-local if it is fibrant, in the sense of the model structure $\bcM $, and for any morphism $g:A\to B$ between cofibrant objects in $S$ the induced morphism from $\bHom (B,Z)$ to $\bHom (A,Z)$ is a weak equivalence in $\SSets $. A morphism $f:X\to Y$ in $\bcC $ is an $S$-local equivalence if the induced morphism from $\bHom (Y,Z)$ to $\bHom (X,Z)$ is a weak equivalence in $\SSets $ for any $S$-local object $Z$ in $\bcC $. Then there exists a new left proper cellular model structure $\bcM _S=(W_S,C_S,F_S)$ on the same category $\bcC $, such that $C_S=C$, $W_S$ consists of $S$-local equivalences in $\bcC $, so contains $W$, and $F_S$ is standardly defined by the right lifting property and so is contained in $F$. The model structure $\bcM _S$ is again left proper and cellular with the same set of generating cofibrations $I$ and the new set of generating trivial cofibrations $J_S$. The model category $(\bcC ,\bcM _S)$ is called the (left) Bousfield localization of $(\bcC ,\bcM _S)$ with respect to $S$. This all can be found in \cite{Hirsch}.

Notice that the identity adjunction on $\bcC $ is a Quillen adjunction and induces the derived adjunction $L\Id :Ho(\bcC )\dashv Ho(\bcC _S):R\Id $, where $Ho(\bcC _S)$ is the homotopy category of $\bcC $ with respect to the model structure $\bcM _S$. Since cofibrations remain the same and, according to our assumption, all objects are cofibrant, the functor $L\Id $ is the identity on objects and surjective on Hom-sets. To describe $R\Id $ we observe the following. Since $F_S$ is smaller than $F$, the fibre replacement functor in $(\bcC ,\bcM )$ is different from the fibre replacement functor in $(\bcC ,\bcM _S)$. Taking into account that $\bcC $ is left proper and cellular, one can show that there exists a fibrant replacement $\Id _{\bcC }\to L_S$ in $\bcM _S$, such that, if $X$ is already fibrant in $\bcM $, then $L_S(X)$ can be more or less visibly constructed from $X$ and $S$, see Section 4.3 in \cite{Hirsch} (or less abstract presentation in \cite{DrorFarjoun}). The right derived functor $R\Id $, being the composition of $Ho(L_S)$ and the functor induced by the embedding of $S$-local, i.e. cofibrant in $\bcM _S$, objects into $\bcC $, identifies $Ho(\bcC _S)$ with the full subcategory in $Ho(\bcC )$ generated by $S$-local objects of $\bcC $.

Since $(\bcC ,\bcM )$ is a simplicial model category, then so is $(\bcC ,\bcM _S)$, see Theorem 4.1.1 (4) in \cite{Hirsch}. Suppose that $\bcC $ is, moreover, closed symmetric monoidal with product $\otimes $, and that the monoidal structure is compatible with the model one in the standard sense, i.e. $\bcC $ is a symmetric monoidal model category (see above). The new model category $(\bcC ,\bcM _S)$ is monoidal model, i.e. the model structure $\bcM _S$ is compatible with the existing monoidal product $\otimes $, if and only if for each $f$ in $S$ and any object $X$ in the union of the domains $\dom (I)$ and codomains $\codom (I)$ of generating cofibrations $I$ in $\bcC $ the product $\id _X\otimes f$ is in $W_S$.


This is exactly the case when the set $S$ is generated by a morphism $p:A\to \uno $, where $A$ is an object in $\bcC $ and $\uno $ is the unit object for the monoidal product $\otimes $, i.e.
  $$
  S=\{ X\wedge A\stackrel{\id _X\wedge p}{\lra }X\; |\; X\in \dom(I)\cup \codom (I)\} \; .
  $$
In that case the model structure $\bcM _S$ is compatible with the monoidal one, so that $(\bcC ,\bcM _S)$ is a simplicial closed symmetric monoidal model category, which is left proper cellular.


Let us write $\bcM _A$ and $L_A$ instead of, respectively, $\bcM _A$ and $L_S$ when $S$ is generated by $A$ in the above sense. One of the fundamental properties of the localization functor $LA$ is that, for any two objects $X$ and $Y$ in $\bcC $, the object $L_A(X\times Y)$ is weak equivalent to the object $L_A(X)\times L_A(Y)$, in the sense of the model structure $\bcM _A$. The proof of this fact in topology is given on page 36 of the book \cite{DrorFarjoun}, and it can be verbally transported to abstract setting. All we need is the Quillen adjunction in two variables in $\bcC $, and the fact saying that if $Y$ is $S$-local, then $\cHom (X,Y)$ is $S$-local for any $X$ in $\bcC $, which is also the consequence of adjunction.

\bigskip

\begin{small}

\end{small}

\bigskip

\bigskip

\begin{small}

{\sc Department of Mathematical Sciences, University of Liverpool, Peach Street, Liverpool L69 7ZL, England, UK}

\end{small}

\medskip

\begin{footnotesize}

{\it E-mail address}: {\tt vladimir.guletskii@liverpool.ac.uk}

\end{footnotesize}

\bigskip

\end{document}